\documentclass[12pt,oneside,reqno]{amsart}
\usepackage{amssymb}
\usepackage{color}
\usepackage{txfonts}
\usepackage{bbm}
\usepackage{cases}
\usepackage{amsmath}
\usepackage{graphicx}
\usepackage{mathrsfs}
\usepackage{stmaryrd}
\usepackage{amsfonts}
\usepackage{enumerate,amsmath,amssymb,amsthm}

\usepackage[nobysame]{amsrefs}
\BibSpec{article}{%
+{}{\PrintAuthors} {author}
+{,}{ \textrm} {title}
+{.}{ \textit} {journal}
+{,}{ \textbf} {volume}
+{}{ \parenthesize} {date}
+{,}{ } {pages}
+{.}{ Available at arXiv:} {eprint}
+{.}{} {transition}
}
\BibSpec{book}{%
+{}{\PrintAuthors} {author}
+{,}{ \textit} {title}
+{.}{ \textrm} {series} %+{,}{ \textrm} {series}
+{,}{ Vol.} {volume} 
+{.}{ } {publisher}
+{,}{ } {date}
+{.}{} {transition}
}
\usepackage{color}
\usepackage[colorlinks=true]{hyperref}
\hypersetup{
    linkcolor=blue,          % color of internal links
    citecolor=red,        % color of links to bibliography
    filecolor=blue,      % color of file links
    urlcolor=cyan  
}

\numberwithin{equation}{section}

\newcommand{\be}{\begin{eqnarray}}
\newcommand{\ee}{\end{eqnarray}}
\newcommand{\ce}{\begin{eqnarray*}}
\newcommand{\de}{\end{eqnarray*}}
\newtheorem{theorem}{Theorem}[section]
\newtheorem{lemma}[theorem]{Lemma}
\newtheorem{remark}[theorem]{Remark}
\newtheorem{definition}[theorem]{Definition}
\newtheorem{proposition}[theorem]{Proposition}
\newtheorem{Examples}[theorem]{Example}
\newtheorem{corollary}[theorem]{Corollary}

\def\nor{{\rm |\mspace{-2mu}[}}
\def\norr{{\rm ]\mspace{-2mu}|}}

\def\eps{\varepsilon}

\def\e{\mathrm{e}}

\def\v{\mathrm{v}}
\def\u{\mathbf{u}}
\def\w{\mathrm{w}}
\def\p{\partial}

\def\[{{\Big[}}
\def\]{{\Big]}}
\def\<{{\langle}}
\def\>{{\rangle}}
\def\({{\Big(}}
\def\){{\Big)}}

\def\bx{{\mathbf{x}}}
\def\tr{\mathrm {tr}}

\def\dif{{\mathord{{\rm d}}}}

\def\min{{\mathord{{\rm min}}}}

\def\no{\nonumber}
\def\={&\!\!=\!\!&}

\def\bB{{\mathbf B}}
\def\bC{{\mathbf C}}

\def\cI{{\mathcal I}}

\def\cL{{\mathcal L}}

\def\cN{{\mathcal N}}

\def\cR{{\mathcal R}}

\def\cT{{\mathcal T}}

\def\mA{{\mathbb A}}
\def\mB{{\mathbb B}}
\def\mC{{\mathbb C}}

\def\mE{{\mathbb E}}
\def\mF{{\mathbb F}}

\def\mI{{\mathbb I}}

\def\mL{{\mathbb L}}

\def\mN{{\mathbb N}}

\def\mP{{\mathbb P}}
\def\mQ{{\mathbb Q}}
\def\mR{{\mathbb R}}

\def\mW{{\mathbb W}}

\def\bB{{\mathbf B}}

\def\1{{\mathbf{1}}}

\def\sA{{\mathscr A}}
\def\sB{{\mathscr B}}
\def\sC{{\mathscr C}}
\def\sD{{\mathscr D}}
\def\sE{{\mathscr E}}
\def\sF{{\mathscr F}}

\def\sI{{\mathscr I}}
\def\sJ{{\mathscr J}}
\def\sK{{\mathscr K}}
\def\sL{{\mathscr L}}

\def\sS{{\mathscr S}}

\def\sU{{\mathscr U}}

\def\geq{\geqslant}
\def\leq{\leqslant}

\def\div{\mathord{{\rm div}}}

\def\eps{\varepsilon}

\def\e{\mathrm{e}}

\def\v{\mathrm{v}}
\def\u{\mathbf{u}}
\def\w{\mathrm{w}}
\def\p{\partial}

\def\[{{\Big[}}
\def\]{{\Big]}}
\def\<{{\langle}}
\def\>{{\rangle}}
\def\({{\Big(}}
\def\){{\Big)}}

\def\bx{{\mathbf{x}}}
\def\tr{\mathrm {tr}}

\def\dif{{\mathord{{\rm d}}}}

\def\min{{\mathord{{\rm min}}}}

\def\no{\nonumber}
\def\={&\!\!=\!\!&}
\def\bt{\begin{theorem}}
\def\et{\end{theorem}}
\def\bl{\begin{lemma}}
\def\el{\end{lemma}}
\def\br{\begin{remark}}
\def\er{\end{remark}}
\def\bx{\begin{Examples}}
\def\ex{\end{Examples}}
\def\bd{\begin{definition}}
\def\ed{\end{definition}}
\def\bp{\begin{proposition}}
\def\ep{\end{proposition}}
\def\bc{\begin{corollary}}
\def\ec{\end{corollary}}

\def\geq{\geqslant}
\def\leq{\leqslant}

\def\div{\mathord{{\rm div}}}

 \def\R{\mathbb R}
 \def\R{\mathbb R}    
  
\def\dd{\delta}   
\def\<{\langle} \def\>{\rangle}

\def\dd{{a}}

\allowdisplaybreaks

\begin{document}

\title{Schauder's estimate for nonlocal kinetic equations\\
 and its applications}
% in SDEs with rough coefficients}

%\date{\today}

\author{Zimo Hao, Mingyan Wu and Xicheng Zhang
}


\thanks{
Research of X. Zhang is partially supported by NNSFC grant of China (No. 11731009)  and the DFG through the CRC 1283 
``Taming uncertainty and profiting from randomness and low regularity in analysis, stochastics and their applications''. }

\address{School of Mathematics and Statistics, Wuhan University, Wuhan, Hubei 430072, P.R.China 
\\ Emails: zimohao@whu.edu.cn, mywu@whu.edu.cn, XichengZhang@gmail.com }

\begin{abstract}
In this paper we develop a new method based on Littlewood-Paley's decomposition 
and heat kernel estimates of integral form, to establish Schauder's estimate for the following degenerate nonlocal equation in $\mathbb R^{2d}$ with H\"older coefficients:
$$
\partial_tu=\mathscr L^{(\alpha)}_{\kappa;{\rm v}} u+b\cdot\nabla u+f,\ u_0=0,
$$
where $u=u(t,x,{\rm v})$ 
and $\mathscr L^{(\alpha)}_{\kappa;{\rm v}}$ is a nonlocal $\alpha$-stable-like operator with $\alpha\in(1,2)$ and kernel function $\kappa$, which acts on the variable ${\rm v}$.
As an application, we show the strong well-posedness to the following degenerate stochastic differential equation with H\"older drift $b$:
$$
{\rm d}Z_t=b(t,Z_t){\rm d}t+(0,\sigma(t,Z_t){\rm d}L^{(\alpha)}_t),\ \ Z_0=(x,{\rm v})\in\mathbb R^{2d},
$$
where $L^{(\alpha)}_t$ is a $d$-dimensional rotationally invariant and symmetric $\alpha$-stable process with $\alpha\in(1,2)$, 
and $b:\mathbb R_+\times\mathbb R^{2d}\to\mathbb R^{2d}$ is a $(\gamma,\beta)$-H\"older continuous function in $(x,{\rm v})$
with $\gamma\in\big(\frac{2+\alpha}{2(1+\alpha)},1\big)$ and $\beta\in\big(1-\frac{\alpha}{2},1\big)$, $\sigma:\mathbb R_+\times\mathbb R^{2d}\to\mathbb R^d\otimes\mathbb R^d$
is a Lipschitz function. 
Moreover, we also show that for almost all $\omega$, the following random transport equation has a unique $C^1_b$-solution:
$$
\partial_tu(t,x,\omega)+(b(t,x)+L^{(\alpha)}_t(\omega))\cdot\nabla_x u(t,x,\omega)=0,\ \ u(0,x)=\varphi(x),
$$
where  $\varphi\in C^1_b(\mathbb R^d)$ and $b:\mathbb R_+\times\mathbb R^d\to\mathbb R^d$ is a bounded continuous function in $(t,x)$ 
and $\gamma$-order H\"older continuous in $x$ uniformly in $t$
with $\gamma\in\big(\frac{2+\alpha}{2(1+\alpha)},1\big)$.
\end{abstract}

\maketitle

\tableofcontents

\section{Introduction}

Let $T>0$. Consider the following backward transport equation (a first order PDE):
\begin{align}\label{TE}
\p_su+b\cdot\nabla_x u+f=0,\ \ u(T,x)=\varphi(x),
\end{align}
where $b(s,x), f(s,x):[0,T]\times\mR^d\to\mR^d$ and $\varphi(x):\mR^d\to\mR$ are measurable functions. 
It is a classical fact that if $b,f,\varphi$ are $C^1_b$-functions in $x$ uniformly with respect to $s$, 
then the above equation is well-posed, and the unique solution is in fact given by
$$
u(s,x)=\varphi(X_{s,T}(x))+\int^T_s f(t,X_{s,t}(x))\dif t,
$$
where for each $x\in\mR^d$, $X_{s,t}(x)$ solves the following ordinary differential equation (abbreviated as ODE):
\begin{align}\label{ODE}
\dot X_{s,t}(x)=b(t,X_{s,t}(x)),\ \ X_{s,s}(x)=x,\ \ t\geq s.
\end{align}
Basing on this representation, DiPerna and Lions \cite{Di-Li} developed a well-posed theory for ODE \eqref{ODE} when
$b\in \mW^{1,1}_{loc}$ (the first order Sobolev space) has bounded divergence through studying the transport 
equation \eqref{TE} (see also  \cite{Am} for the investigation of ODE \eqref{ODE} with BV-vector field $b$). 
The corresponding results for SDEs are referred to \cite{Fi}, \cite{Zh10}, \cite{Fa-Lu-Th} and \cite{Zh13}.
It should be noticed that if $b$ is only H\"older continuous, PDE \eqref{TE} would be ill-posed (see \cite{Fl-Gu-Pr} for counterexamples).
On the other hand, when $b$ is H\"older continuous, under some random perturbations, it was shown in \cite{Fl-Gu-Pr} that the following transport equation (an stochastic PDE) is well-posed:
\begin{align}\label{GA8}
\dif u+(b\cdot\nabla_x u)\dif s+u\circ\dif W_s=0,\ \ u(T,x)=\varphi(x),
\end{align}
where $\circ$ stands for the Stratonovich integral, and $W$ is a standard $d$-dimensional Brownian motion on some probability space.

\medskip

In the same spirit, we consider the following backward  heat equation:
\begin{align}\label{SE}
\p_su+\Delta u+b\cdot\nabla_x u+f=0,\ \ u(T,x)=\varphi(x),
\end{align}
where $\Delta$ is the Laplacian in $\mR^d$.
When $b,f, \varphi$ are $C^2_b$-functions, the unique solution $u$ still has the following representation:
$$
u(s,x)=\mE \varphi(X_{s,T}(x))+\int^T_s \mE f(t,X_{s,t}(x))\dif t,
$$
where for each $x\in\mR^d$, $X_{s,t}(x)$ is the stochastic flow defined by stochastic differential equation (abbreviated as SDE)
\begin{align}\label{SDE}
\dif X_{s,t}(x)=b(s,X_{s,t}(x))\dif s+\sqrt{2}\dif W_s,\ \ X_{s,s}(x)=x,\ \ t\geq s.
\end{align}
Basing on the $L^p$-theory to PDE \eqref{SE},
Krylov and R\"ockner \cite{Kr-Ro} (see also \cite{Zh05}, \cite{Zh11}, \cite{Fe-Fl})
developed a well-posedness theory for SDE \eqref{SDE} with very singular drift $b$, which reveals the regularization effect
of Brownian noises.
In particular, when $b$ is H\"older continuous, it can be shown that $\{X_{s,t}(x), 0\leq s<t<\infty, x\in\mR^d\}$ 
defines a $C^1$-stochastic diffeomorphism flow so that $u(s,x):=\varphi(X_{s,T}(x))$ solves SPDE \eqref{GA8} in generalized sense 
(cf. \cite{Ku}, \cite{Fl-Gu-Pr}).

\medskip

In this work we are concerning with the following degenerate nonlocal equation in $\mR^{2d}$:
\begin{align}\label{GQ1}
\p_tu=\sL^{(\alpha)}_{\kappa;\v} u+b\cdot\nabla u+f,\ u_0=0,
\end{align}
where $u=u(t,x,\v)$ and $\sL^{(\alpha)}_{\kappa;\v}$ is an $\alpha$-stable-like 
operator acting on the variable $\v$ with the form:
\begin{align}\label{KAP}
\sL^{(\alpha)}_{\kappa;\v} u(x,\v):=\int_{\mR^d}(u(x,\v+w)+u(x,\v-w)-2u(x,\v))\frac{\kappa(t,x,\v,w)}{|w|^{d+\alpha}}\dif w,
\end{align}
where $\alpha\in (0,2)$ and $\kappa(t,x,\v,w)$ is symmetric in $w$, and $b(t,x,\v)$ takes the form
\begin{align}\label{GA0}
b(t,x,\v)=(b^{(1)}(t,x,\v),b^{(2)}(t,x,\v)):\mR_+\times\mR^{2d}\to\mR^{2d}.
\end{align}
Notice that the typical example of equation \eqref{GQ1} is the following nonlocal kinetic equation:
\begin{align}\label{Kinetic}
\p_t u=\sK u+f\mbox{ with }\sK u:=\sL^{(\alpha)}_{\kappa;\v}u+\v\cdot\nabla_x u,
\end{align}
 which naturally occurs in the study of spatial inhomeogeneous Boltzmann equations (cf. \cite{Vi}).

\medskip

The first goal of this paper is to establish the following Schauder's type estimate for \eqref{GQ1}:
\begin{align}\label{Scha}
\|u\|_{L^\infty([0,T]; \bC^{(\gamma+\alpha)/(1+\alpha)}_x\cap \bC^{\alpha+\beta}_\v)}\leq C\|f\|_{L^\infty([0,T]; \bC^{\gamma/(1+\alpha)}_x\cap \bC^{\beta}_\v)},
\end{align}
where $\alpha\in(1,2)$ and $\beta\in(0,1)$, $\gamma\in[\beta,1+\alpha)$. Here $b$ and $\kappa$ satisfy some H\"older assumptions 
(see {\bf (H$^{\alpha,\vartheta}_{\beta,\gamma}$)} below for precise statement).
In PDE's theory, Schauder's estimate plays a basic role in constructing the classical solution for quasilinear PDEs. Nowadays, 
there are many ways to prove such an estimate for heat equations 
(see \cite{Gi-Tr}, \cite{Kr}, \cite{Kr-Pr17}). In recent years, 
Schauder's estimate for nonlocal equations are also drawn great interests (see \cite{Ba}, \cite{Ba-Ka}, 
\cite{Do-Ki}, \cite{Im-Ji-Sh}, \cite{Zh-Zh18}, etc.). 
However, most of the works are concentrated on the non-degenerate case.
In the degenerate case, Lunardi \cite{Lu} showed Schauder's estimate in anisotropic H\"older spaces
for linear degenerate Kolmogorov's equations. 
Here it is natural to use the anisotropic H\"older spaces due to the feature of multiple scales in different directions.
Later, in \cite{Lo05} and \cite{Pr09}, the authors established Schauder's estimate for hypoelliptic Kolmogorov equations with partial nonlinear smooth drifts 
(corresponding to $b^{(1)}(t,x,\v)=\v$ in \eqref{GA0}).
For general variable coefficient $b$ case, to the best of our knowledge, 
the authors in \cite{Ch-Ho-Me1} first establish the sharp Schauder estimate for degenerate nonlinear Kolmogorov equations
under some weak H\"ormander's conditions, which  in our case corresponds to \eqref{Scha} with
$\alpha=2$ and $\gamma=\beta$. 
Their method is based on complex forward parametrix expansions.
We mention that the $L^p$-maximal regularity for degenerate nonlocal Kolmogorov's equations with constant coefficients was also obtained 
in \cites{Ch-Zh18, Ch-Zh17}, \cite{Hu-Me-Pr}.

\medskip

To establish Schauder's estimate \eqref{Scha}, we develop a completely new method, that is based on 
Littlewood-Paley's decomposition and heat kernel estimates of some integral forms. Roughly speaking, 
when we consider the usual heat equation, due to Besov's characterization of H\"older spaces, 
the key point is the following integral 
form estimate of  the heat kernel (see Lemma \ref{Le11} below): for any $\beta\geq 0$ and some constant $C=C(d,\beta)>0$,
$$
\int^t_0\left(\int_{\mR^d}|x|^\beta |\cR_j p_{s}(x)|\dif x\right)\dif s\leq C 2^{-2j-\beta j},\ \ \forall t\geq 0,\ \ j\in\mN,
$$
where $\cR_j$ is the usual block operator in Littlewood-Paley's decomposition, and $p_s(x)$ 
is the Gaussian heat kernel. Unlike the usual method by firstly showing Schauder's estimate 
for constant coefficient equations, then freezing it for variable coefficient equations, we directly do it for variable coefficient
equations by Duhamel's formula (see Theorem \ref{Th42} below), which looks simpler. Moreover, the advantage of our method
is that it provides more flexibility to borrow the spatial regularity of coefficients
to compensate the time singularity when we use it to treat the degenerate equation, which allows us to obtain the sharp Schauder estimate \eqref{Scha}.

\medskip

Another goal of this paper is to use \eqref{Scha} to show the strong well-posedness 
as well as the $C^1$-stochastic diffeomorphism flow property to degenerate SDEs driven by $\alpha$-stable processes with H\"older drifts. 
In particular, we shall prove the well-posedness to the following random transport equation with H\"older coefficient:
\begin{align}\label{JQ1}
\p_tu(t,x,\omega)+(b(t,x)+L^{(\alpha)}_t(\omega))\cdot\nabla_x u(t,x,\omega)=0,\ \ u(0,x)=\varphi(x),
\end{align}
where $b:\mR_+\times\mR^d\to\mR^d$ is a bounded continuous function in $(t,x)$ 
and $\gamma$-order H\"older continuous function in $x$ uniformly in $t$
with $\gamma\in\big(\frac{2+\alpha}{2(1+\alpha)},1\big)$ 
and $\varphi\in C^1_b(\mR^d)$. Here $L^{(\alpha)}_t$ is a symmetric and rotationally invariant $\alpha$-stable process with $\alpha\in(1,2)$.
Compared with Flandoli, Gubinelli and Priola's work \cite{Fl-Gu-Pr} for {\it stochastic} PDE \eqref{GA8},
it is a little surprise that as a deterministic equation, PDE \eqref{JQ1} would be ill-posed for each fixed $\omega$, while in the pathwise sense, {\it random} PDE \eqref{JQ1} could be well-posed for almost all $\omega$ (see Theorem \ref{Th75} below).

\medskip

In the nondegenerate Brownian diffusion case, as mentioned above, there are numerous works devoted to the studies of strong and weak well-posedness for the SDEs with singular and even distributional drifts (see \cite{Kr-Ro}, \cite{Zh05},  \cite{Fe-Fl}, \cite{Zh11}, \cite{Zh-Zh18a} and references therein). 
While in the nondegenerate and $\alpha$-stable noise case, 
recently there are also several works (see \cite{Pr12}, \cite{Ch-So-Zh}, \cite{Ch-Zh-Zh}, \cite{Cha-Me-Pr}) 
to study this problem, especially for the supercritical case $\alpha\in(0,1)$,
because in this case, from the view point of PDEs, the drift term plays a dominant role. 
On the other hand, in the degenerate Brownian diffusion case, 
Chaudru \cite{Cha} showed a strong uniqueness result for SDEs with H\"older drifts (see also \cite{Wa-Zh15}, \cite{Wa-Zh16}). 
More recently, Chaudru, Honor\'e and Menozii \cite{Ch-Ho-Me2} applied their Schauder's estimate \cite{Ch-Ho-Me1} 
to establish the strong uniqueness for a chain of oscillators driven by Brownian noises.
However, it seems that there are few works to study the {\it degenerate} SDEs with jumps and H\"older drifts.

\medskip

This paper is organized as follows: In Section 2, we recall the well-known anisotropic Besov and H\"older-Zygmund spaces for later use.
In Section 3, we introduce the basic idea of using Littlewood-Paley's decompostion to establish Schauder's estimate for heat equations with variable coefficients.
In Section 4, we prove several commutator estimates, which plays a crucial role in showing the Schauder estimate \eqref{Scha}.
In Section 5, we give the heat kernel estimate of integral form for nonlocal kinetic operators, which is the basic tool for proving Schauder's estimate.
In Section 6, we prove Schauder's estimate \eqref{Scha} under some natural H\"older's assumptions on $\kappa$ and $b$ (see Theorem \ref{Th66}).
In Section 7, we apply the Schauder estimate to the well-posedness of degenerate SDEs with H\"older drifts and also show the well-posedness of
a random transport equation with H\"older drift. The key point is to establish the $C^1$-stochastic diffeomorphism flow property to the degenerate SDEs.
Finally, in Section 8 we show the existence of smooth solutions for degenerate nonlocal equations with unbounded coefficients by a purely probabilistic argument, which has 
independent interest. Throughout this paper we shall use the following conventions:
\begin{itemize}
\item We use $A\lesssim B$  to denote $A\leq C B$ or some unimportant constant $C>0$. 
\item We use  $A\asymp B$ to denote $C^{-1}B\leq A\leq CB$ for some unimportant constant $C\geq 1$. 
\item For any $\eps\in(0,1)$, we use $A\lesssim\eps B+D$ to denote
$A\leq \eps B+C_\eps D$ for some constant $C_\eps>0$.
\item For two operators $\sA_1, \sA_2$, we use $[\sA_1,\sA_2]:=\sA_1\sA_2-\sA_2\sA_1$ to denote their commutator.
\item For a Banach space $\mB$ and $T>0$, we denote
$$
\mL^\infty_T(\mB):=L^\infty([0,T];\mB),\ \mL^\infty_{loc}(\mB):=\cap_{T>0}\mL^\infty_T(\mB),\ \mL^\infty_T:=L^\infty([0,T]\times\mR^d).
$$
\item $\mN_0:=\mN\cup\{0\}$, $\mR_+:=[0,\infty)$, $a\vee b:=\max(a,b)$, $a\wedge b:=\min(a,b)$.

\end{itemize}

\section{Anisotropic Besov and H\"older-Zygmund spaces}

We first introduce the H\"older (and H\"older-Zygmund) spaces. For $h\in\mR^d$ and $f:\mR^d\to\mR$, the first order difference operator is defined by
$$
\delta_hf(x):=f(x+h)-f(x).
$$
For $\beta>0$, let $\sC^{\beta}$ be the usual $\beta$-order H\"older space consisting of all functions $f:\mR^d\to\mR$ with
$$
\|f\|_{\sC^\beta}:=\|f\|_\infty+\cdots+\|\nabla^{[\beta]} f\|_\infty+[\nabla^{[\beta]}  f]_{\sC^{\beta-[\beta]}}<\infty,
$$
where $[\beta]$ denotes the greatest integer less than $\beta$, and $\nabla^j$ stands for the $j$-order gradient, and
$$
[f]_{\sC^\gamma}:=\sup_h\|\delta_h f\|_\infty/|h|^{\gamma},\ \gamma\in[0,1).
$$ 
The $\beta$-order H\"older-Zygmund space $\bC^\beta$ is defined by
$$
\|f\|_{\bC^\beta}:=\|f\|_\infty+[f]_{\bC^\beta}<\infty,\ [f]_{\bC^\beta}:=\sup_h\|\delta^{[\beta]+1}_hf\|_\infty/|h|^{\beta},
$$
where for an integer $m$, $\delta^m_h:=\delta_h\cdots\delta_h$ denotes the $m$-order difference operator.
Notice that for $0<\beta\notin\mN$ and $m\in\mN$ (cf. \cite{Tr83}), 
\begin{align}\label{DG1}
\|f\|_{\sC^\beta}\asymp\|f\|_{\bC^\beta},\ \|f\|_{\bC^m}\lesssim\|f\|_{\sC^m}.
\end{align}
Let 
$$
\<f,g\>:=\int_{\mR^d} f(x) g(x)\dif x.
$$
The adjoint operator of $\delta_h$ with respect to the above $\<\cdot,\cdot\>$  is given by
$$
\delta^*_h=-\delta_{-h}\Leftrightarrow\<\delta_h f, g\>=\<f,\delta^*_h g\>.
$$
In particular, we have
\begin{align}\label{DU}
\delta^*_h\delta_h f(x)=f(x+h)+f(x-h)-2f(x),
\end{align}
and for any  $f\in \sC^2$,
\begin{align}
\|\delta^*_h\delta_h f\|_\infty&\leq (2\|\nabla^2 f\|_\infty|h|^2)\wedge(4\|f\|_\infty).\label{EV31}
\end{align}

Let $\sS(\mR^d)$ be the Schwartz space of all rapidly decreasing functions on $\mR^d$, and $\sS'(\mR^d)$ 
the dual space of $\sS(\mR^d)$ called Schwartz generalized function (or tempered distribution) space. Given $f\in\sS(\mR^d)$, 
the Fourier transform $\hat f$ and inverse Fourier transform $\check f$ are defined by
\begin{align*}
\hat f(\xi)&:=\frac{1}{\sqrt{2\pi}}\int_{\mR^d} \e^{-i\xi\cdot x}f(x)\dif x, \quad\xi\in\mR^d,\\
\check f(x)&:=\frac{1}{\sqrt{2\pi}}\int_{\mR^d} \e^{i\xi\cdot x}f(\xi)\dif\xi, \quad x\in\mR^d.
\end{align*}
Let $m=(m_1,\cdots,m_n)\in\mN^n$ with $m_1+\cdots+m_n=d$ and $a=(a_1,\cdots,a_n)\in[1,\infty)^n$ be fixed.
We introduce the following distance in $\mR^d$ by
$$
|x-y|_a:=\sum_{i=1}^n|x_i-y_i|^{1/a_i},\ x_i,y_i\in\mR^{m_i}.
$$
For $x=(x_1,\cdots, x_n)$, $t>0$ and $s\in\mR$, we denote
$$
t^{s a} x:=(t^{s a_1}x_1,\cdots, t^{s a_n}x_n)\in\mR^d,\ \ B^a_t:=\Big\{x\in\mR^d: |x|_a\leq t\Big\}.
$$
Clearly we have
$$
|t^a x|_a=t|x|_a,\ \ t>0.
$$
Let $\phi^a_0$ be a radial $C^\infty$-function on $\mR^d$ with 
$$
\phi^a_0(\xi)=1\ \mbox{ for } \ \xi\in B^a_1\ \mbox{ and }\ \phi^a_0(\xi)=0\ \mbox{ for } \ \xi\notin B^a_2.
$$
For $\xi=(\xi_1,\cdots,\xi_n)\in\mR^d$ and $j\in\mN$, define 
$$
\phi^\dd_j(\xi):=\phi^a_0(2^{-a j}\xi)-\phi^a_0(2^{-a(j-1)}\xi).
$$
It is easy to see that for $j\in\mN$, $\phi^a_j(\xi)=\phi^a_1(2^{-a(j-1)}\xi)\geq 0$ and
$$
{\rm supp}\phi^a_j\subset B^a_{2^{j+1}}\setminus B^a_{2^{j-1}},\ \  \sum_{j=0}^{k}\phi^\dd_j(\xi)=\phi^a_0(2^{-a k}\xi)\to 1,\ \ k\to\infty.
$$
\bd[Anisotropic Besov and H\"older-Zygmund spaces]\label{iBesov}
For given $j\in\mN_0$, the block operator $\cR^\dd_j$ is defined on $\sS'(\mR^d)$ by
\begin{align}\label{Def2}
\cR^\dd_jf(x):=(\phi^\dd_j\hat f)\check{\,\,}(x)=\check\phi^\dd_j* f(x)=2^{a\cdot m (j-1)}\int_{\mR^d}\check\phi^a_1(2^{a(j-1)}y) f(x-y)\dif y,
\end{align}
where $a\cdot m=a_1m_1+\cdots+a_nm_n$. 
For any $s\in\mR$, the anisotropic Besov space $\bB^s_{\dd,\infty}$ is defined by
$$
\bB^s_{\dd,\infty}:=\Big\{f\in\sS'(\mR^d): \|f\|_{\bB^s_{\dd,\infty}}:= \sup_{j \geq 0}\left( 2^{s j} \|\cR^\dd_j f\|_{\infty} \right)<\infty\Big\},
$$
%{\rm (Anisotropic H\"older spaces)}For any $s\in(0,2)$, 
and for $s\geq 0$, the anisotropic H\"older-Zygmund space $\bC^s_{\dd}$ is defined by
$$
\bC^s_{\dd}:=\Big\{f\in \mR^d\to\mR: \|f\|_{\bC^s_{\dd}}:= \|f\|_\infty+[f]_{\bC^s_{\dd}}<\infty\Big\},
$$
where 
$$
[f]_{\bC^s_{\dd}}:=\sup_{h}\|\delta^{[s]+1}_h f\|_{\infty}/|h|^s_a.
$$
In particular, if $a=(1,\cdots,1)$, we shall drop the index $\dd$ in $\bB^s_{\dd,\infty}$, $\cR^\dd_j$ and $\bC^s_a$. 
\ed

For $j\in\mN_0$, by definition it is easy to see that
\begin{align}\label{KJ2}
\cR^a_j=\cR^a_j\widetilde\cR^a_j,\ \mbox{ where }\ \widetilde\cR^a_j:=\cR^a_{j-1}+\cR^a_{j}+\cR^a_{j+1}\mbox{ with } \cR^a_{-1}\equiv 0,
\end{align}
and $\cR^a_j$ is symmetric in the sense that
$$
\<\cR^a_j f,g\>=\< f,\cR^a_jg\>.
$$
The cut-off low frequency operator $S_k$ is defined by
\begin{align}\label{EM9}
S_kf:=\sum_{j=0}^{k-1}\cR^a_j f=2^{a\cdot m k}\int_{\mR^d}\check\phi^a_0(2^{ka}(x-y))f(y)\dif y\to f.
\end{align}
For $f,g\in\sS'(\R^d)$, define
$$
T_fg=\sum_{k\geq 2} S_{k-1}f\cR^a_k g,\ \ R(f,g):=\sum_{k\in\mN}\sum_{|i|\leq1}\cR^a_k f\cR^a_{k-i}g.
$$
The Bony decomposition of $fg$ is formally given by (cf. \cite{Ba-Ch-Da})
\begin{align}\label{Bony}
fg=T_fg+T_g f+R(f,g).
\end{align}
The key point of Bony's decomposition is 
\begin{align}\label{YQ1}
\cR^\dd_j (S_{k-1}f\cR^\dd_k g)=0 \ \mbox{ for }\ |k-j|>4.
\end{align}
Indeed, by Fourier's transform, we have
\begin{align*}
\big(\cR^\dd_j (S_{k-1}f\cR^\dd_k g)\big)\,{\hat{}}=\phi^\dd_j\cdot \Big((\phi_0(2^{a(2-k)}\cdot) \hat f)*(\phi^\dd_k\hat g)\Big).
\end{align*}
Since the support of $(\phi_0(2^{a(2-k)}\cdot) \hat f)*(\phi^\dd_k\hat g)$ is contained in $B^a_{2^{k+2}}\setminus B^a_{2^{k-2}}$, if $|k-j|>4$, then
$$
\phi^\dd_j\cdot \Big((\phi_0(2^{a(2-k)}\cdot) \hat f)*(\phi^\dd_k\hat g)\Big)=0.
$$

The following result gives the equivalence between $\bB^s_{\dd,\infty}$ and $\bC^s_a$ (cf. \cite{Tr92}, \cite{Dac}).
\bt\label{Th22}
For any $s>0$, it holds that
\begin{align}\label{LJ1}
\|f\|_{\bB^s_{\dd,\infty}}\asymp\|f\|_{\bC^s_{\dd}}\asymp \|f\|_{\bC^{s/a_1}_{x_1}}+\cdots+\|f\|_{\bC^{s/a_n}_{x_n}},
\end{align}
where $\|f\|_{\bC^{s/a_i}_{x_i}}:=\sup_{x_j\in\mR^{m_j},j\not=i}\|f(x_1,\cdots,x_{i-1},\cdot,x_{i+1},\cdots, x_n)\|_{\bC^{s/a_i}}$.
By convention we denote
$$
\bC^s_\dd:=\bB^s_{a,\infty},\ \ s<0.
$$
\et
We have the following interpolation inequality.
\bc
For any $s<r<t$, there is a constant $C>0$ such that for any $\eps\in(0,1)$,
\begin{align}
\|f\|_{\bC^{r}_a}\leq C\|f\|^{(t-r)/(t-s)}_{\bC^{s}_a}\|f\|^{(r-s)/(t-s)}_{\bC^{t}_a}\leq \eps\|f\|_{\bC^t_a}+C\eps^{(s-r)/(t-r)}\|f\|_{\bC^s_a}.\ \ \label{In}
\end{align}
\ec
\begin{proof}
By \eqref{LJ1} and the definition of $\bB^r_{a,\infty}$, we have
\begin{align*}
\|f\|_{\bC^{r}_a}\lesssim \|f\|_{\bB^{r}_{a,\infty}}=\sup_{j\geq 0}2^{rj}\|\cR^a_j f\|_\infty
&\leq\left(\sup_{j\geq 0}2^{sj}\|\cR^a_j f\|_\infty\right)^{(t-r)/(t-s)}\left(\sup_{j\geq 0}2^{tj}\|\cR^a_j f\|_\infty\right)^{(r-s)/(t-s)}\\
&=\|f\|^{(t-r)/(t-s)}_{\bB^{s}_{a,\infty}}\|f\|^{(r-s)/(t-s)}_{\bB^{t}_{a,\infty}}\lesssim \|f\|^{(t-r)/(t-s)}_{\bC^{s}_a}\|f\|^{(r-s)/(t-s)}_{\bC^{t}_a}.
\end{align*}
The desired interpolation inequality follows.
\end{proof}

\section{Schauder's estimates for heat equations}

In this section we present the basic idea of proving Schauder's estimate for heat equations by Littlewood-Paley's theory.
Let $(a^{ij}(t))$ be a measurable symmetric matrix-valued function on $\mR_+$ and satisfy that for some $c_0\geq 1$,
\begin{align}\label{CN1}
c^{-1}_0|\xi|^2\leq a^{ij}(t)\xi_i\xi_j\leq c_0|\xi|^2,\ \ \xi\in\mR^d,\ t\geq 0. 
\end{align}
Define for $0\leq s<t<\infty$ and $x\in\mR^d$,
\begin{align}\label{X0}
p_{s,t}(x):=\frac{\e^{-\<A_{s,t}^{-1}x,x\>/2}}{\sqrt{2\pi\det (A_{s,t})}}=\frac{\e^{-\<\bar A_{s,t}^{-1}x,x\>/(2(t-s))}}{\sqrt{2\pi(t-s)^d\det (\bar A_{s,t})}},
\end{align}
where 
$$
A_{s,t}:=\int^t_s a(r)\dif r=(t-s)\int^1_0 a(s+(t-s)r)\dif r=:(t-s) \bar A_{s,t}.
$$

The following lemma is the key observation for Schauder's estimate of heat equation.
\bl\label{Le11}
Under \eqref{CN1}, for any $\beta\geq 0$, 
there is a constant $C=C(c_0,\beta,d)>0$ such that for all $t\geq 0$ and $j\in\mN$,
\begin{align}\label{GF02}
\int^t_0\left(\int_{\mR^d}|x|^\beta |\cR_j p_{s,t}(x)|\dif x\right)\dif s\leq C 2^{-2j-\beta j}.
\end{align}
\el
\begin{proof}
We first show that for any $m\in\mN_0$ and $\beta\geq 0$, 
there is a constant $C=C(c_0,\beta,m,d)>0$ such that for all $0\leq s<t<\infty$ and $j\in\mN$,
\begin{align}\label{GF1}
\int_{\mR^d}|x|^\beta |\cR_j p_{s,t}(x)|\dif x\leq C 2^{-2jm}(t-s)^{-m}\Big(2^{-j}+(t-s)^{1/2}\Big)^\beta.
\end{align}
Recalling \eqref{Def2} and by the change of variables, we have
$$
\int_{\mR^d}|x|^\beta |\cR_j p_{s,t}(x)|\dif x=2^{-jd-j\beta}\int_{\mR^d}|x|^\beta\left|\int_{\mR^d}p_{s,t}(2^{-j}(x-y))\check\phi_1(y)\dif y\right|\dif x.
$$
Since the support of $\phi_1$ is contained in the annulus, by Fourier's transform we have
$$
\int_{\mR^d}p_{s,t}(2^{-j}(x-y))\check\phi_1(y)\dif y=\int_{\mR^d}\Delta^{m}p_{s,t}(2^{-j}(x-\cdot))(y)\cdot\Delta^{-m}\check\phi_1(y)\dif y,\ m\in\mN_0,
$$
where  $\Delta^{-m}\check\phi_1:=(|\xi|^{-2m}\phi_1(\xi))\check{\,}$.
Moreover, by \eqref{X0} and elementary calculations, we have
$$
2^{-j\beta-jd}\int_{\mR^d}|x|^\beta|\Delta^m p_{s,t}(2^{-j}\cdot)(x)|\dif x\leq C2^{-2jm}(t-s)^{\beta/2-m}.
$$
Hence,
\begin{align*}
\int_{\mR^d}|x|^\beta |\cR_j p_{s,t}(x)|\dif x&\lesssim 2^{-jd-j\beta}\int_{\mR^d}|x|^\beta|\Delta^m p_{s,t}(2^{-j}x)|\dif x\int_{\mR^d}|\Delta^{-m}\check\phi^\dd_1(y)|\dif y\\
&+2^{-jd-j\beta}\int_{\mR^d}|\Delta^mp_{s,t}(2^{-j}x)|\dif x\int_{\mR^d}|y|^\beta|\Delta^{-m}\check\phi^\dd_1(y)|\dif y\\
&\lesssim 2^{-2jm}(t-s)^{\beta/2-m}+2^{-j\beta-2jm}(t-s)^{-m},
\end{align*}
which in turn gives \eqref{GF1}.

Let $\sI$ be the left hand side of \eqref{GF02}. We make the following decomposition:
\begin{align*}
\sI=\left(\int^t_{t-t\wedge 2^{-2j}}+\int^{t-t\wedge 2^{-2j}}_0\right)\left(\int_{\mR^d}|x|^\beta |\cR_j p_{s,t}(x)|\dif x\right)\dif s=:\sI_1+\sI_2.
\end{align*}
For $\sI_1$, by \eqref{GF1} with $m=0$, we have 
\begin{align*}
\sI_1&\lesssim \int^t_{t-t\wedge 2^{-2j}}\Big(2^{-j}+(t-s)^{-1/2}\Big)^\beta\dif s
=\int^{t\wedge 2^{-2j}}_0\Big(2^{-j}+s^{1/2}\Big)^\beta\dif s\lesssim 2^{-2j-\beta j}.
\end{align*}
For $\sI_2$, by \eqref{GF1} with $m=2$, we have 
\begin{align*}
\sI_2&\lesssim \int^{t-t\wedge 2^{-2j}}_02^{-4j}(t-s)^{-2}\Big(2^{-j}\vee (t-s)^{-1/2}\Big)^\beta\dif s
=2^{-4j} \int^t_{t\wedge 2^{-2j}}s^{-2-\beta/2}\dif s\lesssim 2^{-2j-\beta j}.
\end{align*}
Combining the above two estimates, we obtain \eqref{GF02}.
\end{proof}

Now we consider the following heat equation with variable coefficients:
\begin{align}\label{PDE}
\p_t u=a^{ij}\p_i\p_j u+f,\ \ u(0)=0,
\end{align}
where $a:\mR_+\times\mR^d\to\mR^d\otimes\mR^d$ is a measurable symmetric matrix-valued function and satisfies 
\begin{enumerate}[\bf (H$^\beta_a$)]
\item For some $c_0\geq 1$ and $\beta\in(0,1)$, it holds that for all $t\geq 0$ and $x,y,\xi\in\mR^d$,
$$
c^{-1}_0|\xi|^2\leq a^{ij}(t,x)\xi_i\xi_j\leq c_0|\xi|^2,\quad |a(t,x)-a(t,y)|\leq c_0|x-y|^\beta.
$$
\end{enumerate}
Below we use Lemma \ref{Le11} to establish Schauder's estimate for heat equation \eqref{PDE}.
\bt\label{Th42}
Let $\beta\in(0,1)$. Under {\bf (H$^\beta_a$)}, there is a constant $C=C(c_0,\beta,d)>0$ such that for any $T>0$,
and $u\in \mL^\infty_T(\bB^{2+\beta}_\infty)$  with $\p_t u\in \mL^\infty_T(\bB^{\beta}_\infty)$ solving PDE \eqref{PDE},
$$
\|u\|_{\mL^\infty_T(\bC^{2+\beta})}\leq C\Big(\|f\|_{\mL^\infty_T(\bC^{\beta})}+\|u\|_{\mL^\infty_T}\Big).
$$
\et
\begin{proof}
Fix $x_0\in\mR^d$ and define
$$
u_{x_0}(t,x):=u(t,x+x_0),\ \ \tilde a_{x_0}(t,x):=a(t,x+x_0)-a(t,x_0).
$$
It is easy to see that
$$
\p_t u_{x_0}=a^{ij}(t,x_0)\p_j\p_j u_{x_0}+\tilde a_{x_0}\p_j\p_j u_{x_0}+f_{x_0},\ \ u_{x_0}(0)=0.
$$
Let $p^{x_0}_{s,t}$ be defined by \eqref{X0} in terms of $a(t,x_0)$. For a space-time function $f$, define
$$
P^{x_0}_{s,t}f(s,x):=\int_{\mR^d}p^{x_0}_{s,t}(x-y)f(s,y)\dif y.
$$
By Duhamel's formula we have
$$
u_{x_0}(t,x)=\int^t_0P^{x_0}_{s,t}\tr(\tilde a_{x_0}\cdot\nabla^2 u_{x_0})(s,x)\dif s
+\int^t_0P^{x_0}_{s,t}f_{x_0}(s,x)\dif s=:I_1(t,x)+I_2(t,x).
$$
Below, without loss of generality, we assume $x_0=0$ and drop the subscript and superscript $x_0$.
First of all, for $I_1(t,x)$, by {\bf (H$^\beta_a$)} and Lemma \ref{Le11}, we have
\begin{align*}
|\cR_jI_1(t,0)|&\leq \int^t_0|\cR_j P_{s,t}\tr(\tilde a\cdot\nabla^2 u)(s,0)|\dif s
\lesssim \int^t_0\left(\int_{\mR^d}|x|^\beta|\cR_j p_{s,t}(x)|\dif x\right)\dif s\|\nabla^2 u\|_{\mL^\infty_T}\\
&\lesssim 2^{-2j-\beta j}\|\nabla^2 u\|_{\mL^\infty_T}\leq \eps 2^{-2j-\beta j}\|u\|_{\mL^\infty_T(\bB^{2+\beta}_\infty)}+2^{-2j-\beta j}\|u\|_{\mL^\infty_T},
\end{align*}
where $\eps>0$ and the last inequality is due to the interpolation and Young's inequalities.
For $I_2(t,x)$,  by \eqref{KJ2} and Lemma \ref{Le11} again, we have
\begin{align*}
|\cR_jI_2(t,0)|&\leq \int^t_0|\cR_j P^{x_0}_{s,t}f(s,0)|\dif s=\int^t_0\left|\int_{\mR^d}\cR_j p^{x_0}_{s,t}(y)f(s,y)\dif y\right|\dif s\\
&=\int^t_0\left|\int_{\mR^d}\cR_j \widetilde\cR_j p^{x_0}_{s,t}(y)f(s,y)\dif y\right|\dif s=\int^t_0\left|\int_{\mR^d}\widetilde\cR_j p^{x_0}_{s,t}(y)\cR_j f(s,y)\dif y\right|\dif s\\
&\leq \int^t_0\left(\int_{\mR^d}|\widetilde\cR_j p_{s,t}(x)|\dif x\right)\dif s\|\cR_jf\|_{\mL^\infty_T}
\lesssim 2^{-2j-\beta j}\|f\|_{\mL^\infty_T(\bB^\beta_\infty)}.
\end{align*}
Combining the above estimates, we obtain that for any $\eps\in(0,1)$ and $j\in\mN$,
\begin{align}\label{Es1}
2^{j(2+\beta)}|\cR_ju(t,x_0)|=2^{j(2+\beta)}|\cR_ju_{x_0}(t,0)|\lesssim
\eps\|u\|_{\mL^\infty_T(\bB^{2+\beta}_\infty)}+\|u\|_{\mL^\infty_T}+\|f\|_{\mL^\infty_T(\bB^\beta_\infty)}.
\end{align}
Moreover, for $j=0$, it is easy to see that
$$
|\cR_ju(t,x_0)|\leq \|u\|_{\mL^\infty_T}.
$$
Thus by the definition of Besov space, we arrive at
$$
\|u\|_{\mL^\infty_T(\bB^{2+\beta}_\infty)}=\sup_{t\in[0,T]}\sup_{j\in\mN_0}2^{j(2+\beta)}\|\cR_ju(t,\cdot)\|_\infty
\leq \eps\|u\|_{\mL^\infty_T(\bB^{2+\beta}_\infty)}+C_\eps\|u\|_{\mL^\infty_T}+C\|f\|_{\mL^\infty_T(\bB^\beta_\infty)},
$$
which gives the desired estimate by choosing $\eps=1/2$ and Theorem \ref{Th22}.
\end{proof}

\section{Commutator estimates}

In the sequel, we shall only consider the following case of anisotropic Besov spaces:
$$
n=2,\ m_1=m_2=d,\ \ a=(1+\alpha,1),\ \mbox{ where }\  \alpha\in(0,2).
$$
For $h\in\mR^d$ and $f(x,\v):\mR^{2d}\to\mR$, we introduce
\begin{align*}
\delta_{h;1} f(x,\v):=\delta_h f(\cdot,\v)(x),&\ \ \delta_{h;2} f(x,\v):=\delta_h f(x,\cdot)(\v),\\
\cR^x_j f(x,\v):=\cR_j f(\cdot,\v)(x),&\ \ \cR^\v_j f(x,\v):=\cR_j f(x,\cdot)(\v),
\end{align*}
and for $\beta>0$, 
\begin{align*}
\|f\|_{\bC_x^{\beta}}:= \sup_{\v\in\mR^{d}} \|f(\cdot,\v)\|_{\bC^\beta},\ \ &\ \  \|f\|_{\bC_\v^{\beta}}:= \sup_{x\in\mR^{d}}\|f(x,\cdot)\|_{\bC^\beta},\\
\|f\|_{\bB_{x,\infty}^{\beta}}:= \sup_{\v\in\mR^{d}} \|f(\cdot,\v)\|_{\bB^\beta_\infty},\ \ &\ \  \|f\|_{\bB_{\v,\infty}^{\beta}}:= \sup_{x\in\mR^{d}}\|f(x,\cdot)\|_{\bB^\beta_\infty}.
\end{align*}
Moreover, for $\beta\in(0,1)$, we introduce the following semi-norm for later use:
\begin{align}\label{nor}
\nor f\norr_{\bC^{1+\beta}_a}:=[f]_{\bC^{(1+\beta)/(1+\alpha)}_x}+\|\nabla_\v f\|_{\bC^\beta_\v}.
\end{align}
For $\gamma,\beta\geq 0$, we define the mixed norm
$$
\bC^{\gamma}_x\bC^{\beta}_\v:=\left\{f(x,\v): \|f\|_{\bC^{\gamma}_x\bC^{\beta}_\v}:=\|f\|_{\bC^{\gamma}_x}+\|f\|_{\bC^{\beta}_\v}+
\sup_{h,h'}\|\delta^{[\gamma]+1}_{h;1}\delta^{[\beta]+1}_{h';2}  f\|_\infty/(|h|^{\gamma}|h'|^{\beta})<\infty\right\},
$$
and for $\gamma\in\mR$ and $\beta\geq0$,
$$
\bB^{\gamma}_{x,\infty}\bC^{\beta}_\v:=\left\{f(x,\v): \|f\|_{\bB^{\gamma}_{x,\infty}\bC^{\beta}_\v}:=\sup_{j\in\mN_0}2^{\gamma j}
\|\cR^x_j f\|_{\bC^{\beta}_\v}<\infty\right\}.
$$
In particular, by Theorem \ref{Th22}, we have for $\gamma\in\mR$ and $\beta>0$,
\begin{align}\label{HQ7}
\sup_{j,\ell}2^{\gamma j/(1+\alpha)}2^{\beta\ell}\|\cR^x_j\cR^\dd_\ell f\|_\infty
\asymp \sup_{j}2^{\gamma j/(1+\alpha)}\|\cR^x_j f\|_{\bC^\beta_a}
\asymp \|f\|_{\bB^{(\gamma+\beta)/(1+\alpha)}_{x,\infty}}+\|f\|_{\bB^{\gamma/(1+\alpha)}_{x,\infty}\bC^{\beta}_\v},
\end{align}
and for $\gamma>0$ ans $\beta\geq 0$,
$$
\bB^{\gamma}_{x,\infty}\bC^{\beta}_\v\asymp \bC^{\gamma}_x\bC^{\beta}_\v.
$$

We list some easy properties for later use.
\bl\label{Le16}
\begin{enumerate}[(i)]
\item For any $\theta\in[0,1]$ and $\beta,\gamma\geq 0$, it holds that for some $C=C(\theta,\gamma,\beta)>0$,
\begin{align}\label{GA66}
\|f\|_{\bC^{\theta\gamma}_x\bC^{(1-\theta)\beta}_\v}\leq C\|f\|^\theta_{\bC^{\gamma}_x}\|f\|^{1-\theta}_{\bC^{\beta}_\v}.
\end{align}
\item For all $j\in\mN_0$, it holds that for some $C=C(\alpha)>0$,
\begin{align}\label{GA6}
\|\nabla_x\cR^a_jf\|_\infty\leq C 2^{(1+\alpha)j}\|\cR^a_jf\|_\infty,\ \ \|\nabla_\v\cR^a_jf\|_\infty\leq C 2^{j}\|\cR^a_jf\|_\infty.
\end{align}
\item For any $\beta\in(0,2)$, it holds that  for some $C=C(\alpha,\beta)>0$,
\begin{align}\label{GA2}
\|\cR^a_j f\|_\infty\leq C2^{-\beta j}[f]_{\bC^{\beta}_a},\ j\geq 1.
\end{align}
\item For any $\beta\in(0,1\wedge \alpha)$, it holds that for some $C=C(\alpha,\beta)>0$,
\begin{align}\label{GA1}
\|\nabla_\v f\|_{\bC^\beta_a}\leq C\nor f\norr_{\bC^{1+\beta}_a},
\end{align}
where $\nor f\norr_{\bC^{1+\beta}_a}$ is defined by \eqref{nor}.
\end{enumerate}
\el
\begin{proof}
(i) Notice that
$$
\|\delta^{[\gamma]+1}_{h;1}\delta^{[\beta]+1}_{h';2}  f\|_\infty\lesssim\|f\|_{\bC^\gamma_x}|h|^\gamma,\ \ 
\|\delta^{[\gamma]+1}_{h;1}\delta^{[\beta]+1}_{h';2}  f\|_\infty\lesssim\|f\|_{\bC^\beta_\v}|h'|^\beta.
$$
Hence,
$$
\|\delta^{[\gamma]+1}_{h;1}\delta^{[\beta]+1}_{h';2}  f\|_\infty
\lesssim \|f\|^\theta_{\bC^{\gamma}_x}\|f\|^{1-\theta}_{\bC^\beta_\v}|h|^{\theta\gamma}|h'|^{(1-\theta)\beta}.
$$
From this we obtain the desired estimate \eqref{GA66}. 
\\
\\
(ii) It is a direct consequence of \eqref{Def2}.
\\
\\
(iii) Noticing that for $j\geq 1$,
$$
\cR^a_j f(x)=\int_{\mR^d}\check\phi_j^a(h)f(x+h)\dif h,\ \ \int_{\mR^d}\check\phi_j^a(h)\dif h=\phi_j^a(0)=0,
$$
by \eqref{DU} and $\check\phi_j^a(-h)=\check\phi_j^a(h)$, we have
\begin{align*}
\cR^a_j f(x)=\frac{1}{2}\int_{\mR^d}\check\phi_j^a(h)\delta^*_h\delta_h f(x)\dif h.
\end{align*}
Hence,
\begin{align*}
\|\cR^a_j f\|_\infty\leq\frac{1}{2}\sup_h\|\delta^*_h\delta_h f\|_\infty/|h|_a^\beta\int_{\mR^d}\check\phi_j^a(h)|h|_a^\beta\dif h\lesssim 2^{-\beta j}[f]_{\bC^\beta_a}.
\end{align*}
(iv) By Theorem \ref{Th22} and definition, we have
\begin{align*}
\|\nabla_\v f\|_{\bC^\beta_a}&\lesssim \|\nabla_\v f\|_\infty+\sup_{j\in\mN} 2^{\beta j}\|\cR^a_j\nabla_\v f\|_\infty
\stackrel{\eqref{GA6}}{\lesssim}\|\nabla_\v f\|_\infty+\sup_{j\in\mN} 2^{(1+\beta) j}\|\cR^a_j f\|_\infty\\
&\stackrel{\eqref{GA2}}{\lesssim} \|\nabla_\v f\|_\infty+[f]_{\bC^{1+\beta}_a}\leq \|\nabla_\v f\|_\infty+[f]_{\bC^{(1+\beta)/(1+\alpha)}_x}+[f]_{\bC^{1+\beta}_\v}
\lesssim \nor f\norr_{\bC^{1+\beta}_a}.
\end{align*}
The proof is complete.
\end{proof}

We now show several commutator estimates, which are extensions of \cite[Lemma 2.3]{Ch-Zh-Zh}, and will play a key role in showing the
Schauder estimate below.
\bl\label{Le12}
(i) For any $\beta\in(0,\alpha)$ and $\gamma\in(-1-\beta,0]$, there is a constant $C=C(d,\beta,\gamma)>0$ such that 
for all $x,\v\in\mR^d$ and $j\geq 5$,
\begin{align}\label{HQ2}
|[\cR^\dd_j, \tilde f]g|(x,\v)\leq C2^{-j(\gamma+1)}
\Big(2^{-j\beta}+|x|^{\frac{\beta}{1+\alpha}}+|\v|^{\beta}\Big)\nor f\norr_{\bC^{1+\beta}_a}\|g\|_{\bC^{\gamma}_{\dd}},
\end{align}
where $\nor f\norr_{\bC^{1+\beta}_a}$ is defined by \eqref{nor}, and 
$$
\tilde f(x,\v):=f(x,\v)-f(0,0)-\v\cdot\nabla_\v f(0,0).
$$
(ii) For any $\beta\in(0,1)$ and $\gamma\in(-\beta,0]$, there is a constant $C=C(d,\beta,\gamma)>0$ such that
\begin{align}
%\|[\cR^\dd_j, f]g\|_\infty&\leq C2^{-j(\beta+\gamma)}[f]_{\bC^{\beta}_{\dd}}\|g\|_{\bC^{\gamma}_{\dd}},\ j\geq 5,\label{GS2}\\
\|[\cR^x_j, f]g\|_\infty&\leq C2^{-j(\beta+\gamma)}[f]_{\bC^{\beta}_x}\|g\|_{\bC^{\gamma}_{x}},\ j\geq 5.\label{GS1}
\end{align}
\el
\begin{proof}
We only prove (i) since (ii) is similar and easier.
First of all, we have
\begin{align*}
|\tilde f(\bar x,\bar\v)-\tilde f(x,\v)|&=|f(\bar x,\bar\v)-f(x,\v)+(\v-\bar\v)\cdot\nabla_\v f(0,0)|\\
&\leq |x-\bar x|^{\frac{1+\beta}{1+\alpha}}[f]_{\bC^{(1+\beta)/(1+\alpha)}_{x}}
+|\v-\bar \v|\big(|x|^{\frac{\beta}{1+\alpha}}+|\v|^{\beta}+|\bar\v|^{\beta}\big)[\nabla_\v f]_{\bC^{\beta}_\dd}\\
&\lesssim \Big(|x-\bar x|^{\frac{1+\beta}{1+\alpha}}+|\v-\bar \v|\big(|x|^{\frac{\beta}{1+\alpha}}+|\v|^{\beta}+|\bar\v|^{\beta}\big)\Big)\nor f\norr_{\bC^{1+\beta}_a},
\end{align*}
where the last step is due to \eqref{GA1}.
Since for $\gamma=0$, \eqref{HQ2} is easily derived from the above estimate, we assume $\gamma<0$ below.
Noting that by \eqref{EM9},
$$
S_{k-1}\tilde f(x,\v)=2^{(2+\alpha)kd}\int_{\mR^{2d}}\check\phi^a_0(2^{(1+\alpha)k}x',2^k\v')\tilde f(x-x',\v-\v')\dif x'\dif\v',
$$
we have 
\begin{align}\label{HQ1}
|S_{k-1}\tilde f(\bar x,\bar\v)-S_{k-1}\tilde f(x,\v)|
\lesssim  \Big(|x-\bar x|^{\frac{1+\beta}{1+\alpha}}+|\v-\bar \v|\big(2^{-k\beta}+|x|^{\frac{\beta}{1+\alpha}}+|\v|^{\beta}+|\bar\v|^{\beta}\big)\Big)\nor f\norr_{\bC^{1+\beta}_a}.
\end{align}
On the other hand, noting that
\begin{align*}
[\cR^\dd_j, S_{k-1}\tilde f]g(x,\v)=\int_{\mR^{2d}}\check\phi^\dd_j(x-\bar x,\v-\bar \v)
\Big(S_{k-1}\tilde f(\bar x,\bar \v)-S_{k-1}\tilde f(x,\v)\Big)g(\bar x,\bar\v)\dif \bar x\dif\bar\v,
\end{align*}
we have by \eqref{HQ1} and \eqref{Def2},
\begin{align}\label{EE5}
|[\cR^\dd_j, S_{k-1}\tilde f]g|(x,\v)\lesssim \Big(2^{-j(1+\beta)}+2^{-j}\big(2^{-k\beta}+|x|^{\frac{\beta}{1+\alpha}}+|\v|^{\beta}\big)\Big)\nor f\norr_{\bC^{1+\beta}_a}\|g\|_\infty.
\end{align}
Now by using Bony's decomposition \eqref{Bony}, we can write
\begin{align*}
[\cR^\dd_j, \tilde f]g
=[\cR^\dd_j,T_{\tilde f}] g+\cR^\dd_j(T_g \tilde f)-T_{\cR^\dd_jg} \tilde f+\cR^\dd_jR(\tilde f,g)-R(\tilde f,\cR^\dd_jg).
\end{align*}
For the first term, by \eqref{YQ1} we have
\begin{align*}
|[\cR^\dd_j,T_{\tilde f}] g|(x,\v)&=\Bigg|\sum_{|k-j|\leq 4}\Big(\cR^\dd_j (S_{k-1}\tilde f\cR^\dd_k g)
-S_{k-1}\tilde f\cR^\dd_j \cR^\dd_k g\Big)(x,\v)\Bigg|\leq\sum_{|k-j|\leq 4}\big|[\cR^\dd_j, S_{k-1}\tilde f]\cR^\dd_k g\big|(x,\v)\\
&\stackrel{\eqref{EE5}}{\lesssim}\sum_{|k-j|\leq 4}\Big(2^{-j(1+\beta)}+2^{-j}\big(2^{-k\beta}+|x|^{\frac{\beta}{1+\alpha}}+|\v|^{\beta}\big)\Big) \nor f \norr_{\bC^{1+\beta}_{\dd}}
\|\cR^\dd_kg\|_\infty\\
&\quad\lesssim 2^{-j\gamma-j}\Big(2^{-j\beta}+|x|^{\frac{\beta}{1+\alpha}}+|\v|^{\beta}\Big)
\nor f\norr_{\bC^{1+\beta}_a}\|g\|_{\bB^{\gamma}_{\dd,\infty}}.
\end{align*}
%Below without loss of generality we assume $j\geq 5$. 
On the other hand,  we also have
\begin{align*}
\|\cR^\dd_j(T_g \tilde f)\|_\infty
&\stackrel{\eqref{YQ1}}{=}\Bigg\|\sum_{|k-j|\leq 4}\cR^\dd_j(S_{k-1} g\cR^\dd_k \tilde f)\Bigg\|_\infty
\leq \sum_{|k-j|\leq 4}\|\cR^\dd_j(S_{k-1} g\cR^\dd_k\tilde f)\|_\infty\\
&\leq \sum_{|k-j|\leq 4}\|S_{k-1} g\cR^\dd_k\tilde  f\|_\infty\lesssim\sum_{|k-j|\leq 4}\sum_{m\leq k-2}\|\cR^\dd_m g\|_\infty\|\cR^\dd_k\tilde f\|_\infty.
\end{align*}
Since $j\geq 5$ and $|k-j|\leq 4$, by \eqref{GA2}, we further have
\begin{align*}
\|\cR^\dd_j(T_g \tilde f)\|_\infty&\lesssim \|g\|_{\bB^{\gamma}_{\dd,\infty}}[\tilde f]_{\bC^{1+\beta}_a}\sum_{|k-j|\leq 4}\sum_{m\leq k-2}2^{-m\gamma}2^{-k(1+\beta)}
\lesssim \|g\|_{\bB^{\gamma}_{\dd,\infty}}[f]_{\bC^{1+\beta}_a}2^{-j(\gamma+1+\beta)},
\end{align*}
where the last step is due to $\gamma<0$ and $[\tilde f]_{\bC^{1+\beta}_a}\lesssim [f]_{\bC^{1+\beta}_a}$. Similarly, 
\begin{align*}
\|T_{\cR^\dd_jg} \tilde f\|_\infty
&\leq\sum_{k\geq j-2}\|S_{k-1}\cR^\dd_jg\cR^\dd_k \tilde f\|_\infty\leq\sum_{k\geq j-2}\|S_{k-1}\cR^\dd_jg\|_\infty\|\cR^\dd_k\tilde  f\|_\infty\\
&\leq\sum_{k\geq j-2}2^{-k(1+\beta)}[\tilde  f]_{\bC^{1+\beta}_a}\|\cR^\dd_jg\|_\infty
\lesssim 2^{-j(\beta+1+\gamma)}[f]_{\bC^{1+\beta}_a}\|g\|_{\bB^{\gamma}_{\dd,\infty}}.
\end{align*}
Finally, since $1+\beta+\gamma>0$, we have
\begin{align*}
\|\cR^\dd_jR(\tilde f,g)\|_\infty&\stackrel{\eqref{YQ1}}{=}\Bigg\|\sum_{|i|\leq 1, k\geq j-4}\cR^\dd_j(\cR^\dd_k\tilde  f\cR^\dd_{k-i} g)\Bigg\|_\infty
\lesssim \sum_{|i|\leq 1, k\geq j-4}\|\cR^\dd_k\tilde f\|_\infty\|\cR^\dd_{k-i} g\|_\infty\\
&\lesssim \sum_{k\geq j-4}2^{-k(1+\beta+\gamma)}[\tilde f]_{\bC^{1+\beta}_a}\|g\|_{\bB^{\gamma}_{\dd,\infty}}
\lesssim 2^{-j(1+\beta+\gamma)}[f]_{\bC^{1+\beta}_a}\|g\|_{\bB^{\gamma}_{\dd,\infty}},
\end{align*}
and
\begin{align*}
&\|R(\tilde f,\cR^\dd_jg)\|_\infty=\Bigg\|\sum_{|i|\leq 1, |k-j|\leq 1}\cR^\dd_{k-i} \tilde f\cR^\dd_{k}\cR^\dd_j g\Bigg\|_\infty
\lesssim [f]_{\bC^{1+\beta}_a}\|g\|_{\bB^{\gamma}_{\dd,\infty}}2^{-j(1+\beta+\gamma)}.
\end{align*}
Combining the above calculations, we complete the proof. 
\end{proof}

\bl
For any $0<\beta\leq \gamma<1$ and $\eta\in(-\gamma,0]$, there is a  $C=C(\beta,\gamma,\eta)>0$ such that
\begin{align}\label{GP1}
\|[\cR^x_j,f]g\|_{\bC^{\beta}_x}\leq C 2^{(\beta-\gamma-\eta)j}
[f]_{\bC^{\gamma}_x}\|g\|_{\bC^{\eta}_{x}},\ \ j\geq 5.
\end{align}
\el
\begin{proof}
If $\ell\leq j+1$, then by \eqref{GS1},
$$
\|\cR^x_\ell [\cR^x_j,f]g\|_\infty\leq\|[\cR^x_j,f]g\|_\infty\lesssim 2^{-(\gamma+\eta)j} [f]_{\bC^\gamma_x}\|g\|_{\bC^{\eta}_{x}}
\lesssim 2^{-\beta\ell}2^{(\beta-\gamma-\eta)j} [f]_{\bC^\gamma_x}\|g\|_{\bC^{\eta}_{x}}.
$$
If $\ell>j+1$, since $\widehat{\cR^x_\ell \cR^x_j f}=\phi_\ell\phi_j\hat f\equiv0$, we have
$$
\cR^x_\ell [\cR^x_j,f]g=\cR^x_\ell \cR^x_j(fg)-\cR^x_\ell(f\cR^x_j g)
=f\cR^x_\ell \cR^x_jg-\cR^x_\ell(f\cR^x_j g)=-[\cR^x_\ell, f]\cR^x_j g.
$$
Thus by \eqref{GS1} again, we have
$$
\|\cR^x_\ell [\cR^x_j,f]g\|_\infty=\|[\cR^x_\ell, f]\cR^x_j g\|_\infty\leq 2^{-\ell\gamma}[f]_{\bC^{\gamma}_x}\|\cR_j^xg\|_{\infty}
\lesssim 2^{-\ell \beta}2^{(\beta-\gamma-\eta)j}[f]_{\bC^{\gamma}_x}\|g\|_{\bC^{\eta}_{x}}.
$$
Hence,
$$
\|[\cR^x_j,f]g\|_{\bC^{\beta}_x}\lesssim \sup_{\ell\in\mN_0}2^{\beta\ell}\|\cR^x_\ell [\cR^x_j,f]g\|_\infty
\lesssim 2^{(\beta-\gamma-\eta)j} [f]_{\bC^\gamma_x}\|g\|_{\bC^{\eta}_{x}}.
$$
The proof is complete.
\end{proof}

\bl\label{Le33}
Let $\beta,\gamma_2,\theta\in(0,1]$ and $\gamma,\gamma_1\in(0,1+\alpha)$. Under the conditions
\begin{align}\label{GF3}
\gamma\vee\gamma_2<\gamma_1,\ \ \theta\gamma_2<\gamma+\beta\leq (1-\theta)\gamma_1+\theta\gamma_2,\ \ \beta\leq\theta\gamma_2,
\end{align}
there is a constant $C>0$ such that for all $j\geq 5$,
\begin{align}\label{FT1}
\|[\cR^x_j,f] g\|_{\bC^{\beta}_a}\leq 2^{-\frac{\gamma}{1+\alpha}j}
\Big([f]_{\bC^{\gamma_1/(1+\alpha)}_x}+[f]_{\bC^{\gamma_2}_\v}\Big)
\Big(\|g\|_{\bC^{(\gamma+\beta-(1-\theta)\gamma_1-\theta\gamma_2)/(1+\alpha)}_{x}}
+\|g\|_{\bC^{(\gamma-\gamma_1)/(1+\alpha)}_{x}\bC^{\beta}_\v}\Big).
\end{align}
\el
\begin{proof}
First of all, by applying \eqref{GP1} with $(\frac{\beta}{1+\alpha},\frac{\gamma_1}{1+\alpha},\frac{\beta+\gamma-\gamma_1}{1+\alpha})$ 
in place of $(\beta,\gamma,\eta)$, we have
$$
\|[\cR^x_j,f]g\|_{\bC^{\beta/(1+\alpha)}_x}\lesssim 2^{-\frac{\gamma}{1+\alpha}j}
[f]_{\bC^{\gamma_1/(1+\alpha)}_x}\|g\|_{\bC^{(\beta+\gamma-\gamma_1)/(1+\alpha)}_{x}}.
$$
Thus, by definition it suffices to prove
\begin{align}\label{AP1}
\|[\cR^x_j,f]g\|_{\bC^{\beta}_\v}\lesssim \sup_{\ell\in\mN_0}2^{-\ell\beta}\|\cR^\v_\ell[\cR^x_j,f]g\|_\infty\lesssim\mbox{RHS of \eqref{FT1}}.
\end{align}
{\bf (Case: $\ell\leq\frac{j}{1+\alpha}$)}. Since $\gamma_2<\gamma_1$ and $\gamma+\beta<(1-\theta)\gamma_1+\theta\gamma_2\leq\gamma_1$, 
by \eqref{GS1}, we have
\begin{align*}
\|\cR^\v_\ell [\cR^x_j,f]g\|_\infty&\lesssim\|[\cR^x_j,f]g\|_\infty\lesssim 2^{-\frac{\gamma+\beta}{1+\alpha}j} 
[f]_{\bC^{\gamma_1/(1+\alpha)}_x}\|g\|_{\bC^{(\gamma+\beta-\gamma_1)/(1+\alpha)}_{x}}\\
&\lesssim 2^{-\ell\beta}2^{-\frac{\gamma}{1+\alpha}j}  [f]_{\bC^{\gamma_1/(1+\alpha)}_x}
\|g\|_{\bC^{(\gamma+\beta-(1-\theta)\gamma_1-\theta\gamma_2)/(1+\alpha)}_{x}}.
\end{align*}
{\bf (Case: $\ell>\frac{j}{1+\alpha}$)}. Notice that
$$
\cR^\v_\ell [\cR^x_j,f]g=[\cR^x_j,f]\cR^\v_\ell g+[\cR^\v_\ell, [\cR^x_j,f]]g=:I_1+I_2.
$$
For $I_1$, since $\gamma<\gamma_1$, by \eqref{GS1}, we have
\begin{align*}
|I_1|&\lesssim  2^{-\frac{\gamma}{1+\alpha}j}[f]_{\bC^{\gamma_1/(1+\alpha)}_x}\|\cR^\v_\ell g\|_{\bC^{(\gamma-\gamma_1)/(1+\alpha)}_{x}}
\lesssim 2^{-\frac{\gamma}{1+\alpha}j} 2^{-\beta \ell} [f]_{\bC^{\gamma_1/(1+\alpha)}_x}\|g\|_{\bC^{(\gamma-\gamma_1)/(1+\alpha)}_{x}\bC^{\beta}_{\v}}.
\end{align*}
For $I_2$, 
by definition, \eqref{GF3} and \eqref{GS1}, we have
\begin{align*}
|I_2|&=\left|\int_{\mR^{d}}\check\phi_\ell(\bar\v)[\cR^x_j,\delta_{\bar\v;2} f(\cdot,\v)]g(\cdot,\v-\bar\v)\dif \bar\v\right|\\
&\lesssim 2^{-\frac{\gamma+\beta-\theta\gamma_2}{1+\alpha}j}
\int_{\mR^{d}}\check\phi_\ell(\bar\v)[\delta_{\bar\v;2} f]_{\bC^{(1-\theta)\gamma_1/(1+\alpha)}_x}
\|g\|_{\bC^{(\gamma+\beta-(1-\theta)\gamma_1-\theta\gamma_2)/(1+\alpha)}_{x}}\dif \bar\v\\
&\lesssim 2^{-\frac{\gamma+\beta-\theta\gamma_2}{1+\alpha}j}
\left(\int_{\mR^{d}}\check\phi_\ell(\bar\v)|\bar\v|^{\theta\gamma_2}\dif \bar\v\right)[f]_{\bC^{(1-\theta)\gamma_1/(1+\alpha)}_x\bC^{\theta\gamma_2}_\v}
\|g\|_{\bC^{(\gamma+\beta-(1-\theta)\gamma_1-\theta\gamma_2)/(1+\alpha)}_{x}}\\
&\lesssim 2^{-\frac{\gamma+\beta-\theta\gamma_2}{1+\alpha}j} 2^{-\theta\gamma_2\ell}
[f]_{\bC^{(1-\theta)\gamma_1/(1+\alpha)}_x\bC^{\theta\gamma_2}_\v}\|g\|_{\bC^{(\gamma+\beta-(1-\theta)\gamma_1-\theta\gamma_2)/(1+\alpha)}_{x}}\\
&\lesssim 2^{-\frac{\gamma}{1+\alpha}j} 2^{-\beta\ell}\Big([f]_{\bC^{\gamma_1/(1+\alpha)}_x}+[f]_{\bC^{\gamma_2}_\v}\Big)
\|g\|_{\bC^{(\gamma+\beta-(1-\theta)\gamma_1-\theta\gamma_2)/(1+\alpha)}_{x}}.
\end{align*}
Hence, for all $\ell\in\mN_0$ and $j\geq 5$,
$$
\|\cR^\v_\ell [\cR^x_j,f]g\|_\infty\lesssim 2^{-\frac{\gamma}{1+\alpha}j} 2^{-\beta\ell}\Big([f]_{\bC^{\gamma_1/(1+\alpha)}_x}+[f]_{\bC^{\gamma_2}_\v}\Big)
\Big(\|g\|_{\bC^{(\gamma+\beta-(1-\theta)\gamma_1-\theta\gamma_2)/(1+\alpha)}_{x}}+\|g\|_{\bC^{(\gamma-\gamma_1)/(1+\alpha)}_{x}\bC^{\beta}_\v}\Big),
$$
which gives \eqref{AP1}.
The proof is complete.
\end{proof}

\bc\label{Cor34}
Let $\vartheta\in(0,\alpha-1)$ and $0<\beta<\gamma<1+\vartheta$.
For any $\eps\in(0,1)$, there are $\theta>0$ close to zero and constants $C_\eps, C>0$ such that for all $j\geq 5$,
\begin{align}
&\|[\cR^x_j,b\cdot\nabla_x] u\|_{\bC^{\theta\beta}_a}\leq 2^{-\frac{(1-\theta)\gamma}{1+\alpha}j}
\nor b\norr_{\bC^{1+\vartheta}_\dd}\Big(\eps\|u\|_{\bC^{(\alpha+(1-\theta)\gamma+\theta\beta)/(1+\alpha)}_x}
+C_\eps\|u\|_{\bC^{\alpha+\beta}_\v}\Big),\label{JH1}\\
&\qquad\|[\cR^x_j,b\cdot\nabla_\v] u\|_{\bC^{\theta\beta}_a}\leq C 2^{-\frac{(1-\theta)\gamma}{1+\alpha}j}
\Big([b]_{\bC^{\gamma/(1+\alpha)}_x}+[b]_{\bC^\beta_\v}\Big)\|u\|_{\bC^{\alpha+\beta}_\v},\label{JH2}\\
&\qquad\,\,\|[\cR^x_j,\sL^{(\alpha)}_{\kappa;\v}]u\|_{\bC^{\theta\beta}_{\dd}}\leq C2^{-\frac{(1-\theta)\gamma}{1+\alpha}j}
\Big([\kappa]_{\bC^{\gamma/(1+\alpha)}_x}+[\kappa]_{\bC^\beta_\v}\Big)\|u\|_{\bC^{\alpha+\beta}_\v},\label{JH3}
\end{align}
where $\sL^{(\alpha)}_{\kappa;\v}$ is defined by \eqref{KAP}.
\ec
\begin{proof}
Let $\theta\in(0,1)$ be fixed, which will be determined below.
\\
\\
(i) By applying Lemma \ref{Le33} with $(\theta\beta,(1-\theta)\gamma,1+\vartheta,1,\theta\beta)$ in place of $(\beta,\gamma,\gamma_1,\gamma_2,\theta)$, we have
\begin{align}
\|[\cR^x_j,b\cdot\nabla_x] u\|_{\bC^{\theta\beta}_a}
&\lesssim 2^{-\frac{(1-\theta)\gamma}{1+\alpha}j}\nor b\norr_{\bC^{1+\vartheta}_\dd}
\Big(\|\nabla_x u\|_{\bC^{((1-\theta)\gamma-(1-\theta\beta)(1+\vartheta))/(1+\alpha)}_{x}}
+\|\nabla_x u\|_{\bC^{((1-\theta)\gamma-1-\vartheta)/(1+\alpha)}_{x}\bC^{\beta}_\v}\Big)\no\\
&\lesssim 2^{-\frac{(1-\theta)\gamma}{1+\alpha}j}\nor b\norr_{\bC^{1+\vartheta}_\dd}
\Big(\|u\|_{\bC^{(\alpha+(1-\theta)\gamma+\theta\beta-(1-\theta\beta)\vartheta)/(1+\alpha)}_{x}}
+\|u\|_{\bC^{(\alpha+(1-\theta)\gamma-\vartheta)/(1+\alpha)}_{x}\bC^{\beta}_\v}\Big).\label{JH4}
\end{align}
Choosing $\theta>0$ small enough so that
$$
\tfrac{\vartheta(\alpha+\beta)}{\alpha+(1-\theta)\gamma+\theta\beta}>\theta\beta,
$$
by \eqref{GA66} and Young's inequality, for any $\eps\in(0,1)$, there is a constant $C_\eps>0$ such that
$$
\|u\|_{\bC^{(\alpha+(1-\theta)\gamma-\vartheta)/(1+\alpha)}_{x}\bC^{\beta}_\v}\leq
\eps\|u\|_{\bC^{(\alpha+(1-\theta)\gamma+\theta\beta)/(1+\alpha)}_x}+C_\eps\|u\|_{\bC^{\alpha+\beta}_\v},
$$
and also,
$$
\|u\|_{\bC^{(\alpha+(1-\theta)\gamma+\theta\beta-(1-\theta\beta)\vartheta)/(1+\alpha)}_{x}}
\leq \eps\|u\|_{\bC^{(\alpha+(1-\theta)\gamma+\theta\beta)/(1+\alpha)}_x}+\|u\|_\infty.
$$
Substituting these two estimates into \eqref{JH4}, we obtain \eqref{JH1}.
\\
\\
(ii) By Lemma \ref{Le33} with $(\theta\beta,(1-\theta)\gamma,\gamma,\beta,\theta)$ in place of 
$(\beta,\gamma,\gamma_1,\gamma_2,\theta)$, we obtain \eqref{JH2}.
\\
\\
(iii) Recalling \eqref{KAP} and \eqref{DU}, and noticing that
$$
[\cR^x_j,\sL^{(\alpha)}_{\kappa;\v}]u(x,\v)=\int_{\mR^d}
\Big([\cR^x_j, \kappa(\cdot,\cdot,w)]\delta^*_{w;2}\delta_{w;2} u\Big)(x,\v)\frac{\dif w}{|w|^{d+\alpha}},
$$
by Lemma \ref{Le33} with $(\theta\beta,(1-\theta)\gamma,\gamma,\beta,\theta)$ in place of 
$(\beta,\gamma,\gamma_1,\gamma_2,\theta)$, we obtain 
\begin{align*}
&\|[\cR^x_j,\sL^{(\alpha)}_{\kappa;\v}]u\|_{\bB^{\theta\beta}_{\dd,\infty}}
\lesssim \int_{\mR^d}\big\|[\cR^x_j,\kappa(\cdot,\cdot,w)]\delta^*_{w;2}\delta_{w;2} u\big\|_{\bC^{\theta\beta}_a}\frac{\dif w}{|w|^{d+\alpha}}\\
&\qquad\lesssim 2^{-\frac{(1-\theta)\gamma}{1+\alpha}j}\Big([\kappa]_{\bC^{\gamma/(1+\alpha)}_x}+[\kappa]_{\bC^\beta_\v}\Big)
\int_{\mR^d}\|\delta^*_{w;2}\delta_{w;2} u\|_{\bC^{\theta\beta}_\v}\frac{\dif w}{|w|^{d+\alpha}}\\
&\qquad\lesssim 2^{-\frac{(1-\theta)\gamma}{1+\alpha}j}\Big([\kappa]_{\bC^{\gamma/(1+\alpha)}_x}+[\kappa]_{\bC^\beta_\v}\Big)
\|u\|_{\bC^{\alpha+\theta\beta+\eps}_\v}\int_{\mR^d}\frac{(1\wedge|w|^{\alpha+\eps})\dif w}{|w|^{d+\alpha}},
\end{align*}
where $\eps\in(0,(1-\theta)\beta)$,
which in turn yields \eqref{JH3} since $\|u\|_{\bC^{\alpha+\theta\beta+\eps}_\v}\lesssim \|u\|_{\bC^{\alpha+\beta}_\v}$.
\end{proof}

\section{Heat kernel estimates of nonlocal kinetic operators}

In this section we consider the following nonlocal kinetic equation with constant coefficients:
$$
\p_t u=\sL^{(\alpha)}_{\kappa;\v} u+U_t\v\cdot\nabla_x u+f=:\sK u+f,\ \ u(0)=0,
$$
where $\kappa(t,w)$ and $U_t$ are measurable functions and satisfy the following assumptions:
\begin{align}\label{CN0}
c_0^{-1}\leq\kappa(t,w)\leq c_0,\ \ c_0\geq 1,
\end{align}
and
\begin{align}\label{CN00}
c_1:=\|U\|_\infty+\sup_{s<t}\Big((t-s)\|\Pi_{s,t}^{-1}\|\Big)<\infty,\mbox{ where }\ \ \Pi_{s,t}:=\int^t_s U_r\dif r.
\end{align}
It is well known that under \eqref{CN0} and \eqref{CN00}, there is a fundamental solution or heat kernel $p_{s,t}(x,\v)$ to kinetic operator $\p_t-\sK$ so that (see \cite[Lemma 2.5]{Ch-Zh18})
\begin{align}\label{DA3}
u(t,x,\v)=\int^t_0P_{s,t}f(s,x,\v)\dif s:=\int^t_0\Big(\Gamma_{\!s,t}p_{s,t}*\Gamma_{\!s,t}f\Big)(s,x,\v)\dif s,
\end{align}
where operator $\Gamma_{\!s,t}$ is defined by
\begin{align}\label{DA4}
\Gamma_{\!s,t}f(x,\v):=f(x+\Pi_{s,t}\v,\v).
\end{align}
Moreover, for any $\beta,\gamma\geq 0$ with $\beta+\gamma<\alpha$ and $n,m\in\mN_0$, there is a constant  $C>0$ such that
\begin{align}
\int_{\mR^{2d}}|x|^{\beta}|\v|^{\gamma}|\nabla^{n}_{x}\nabla^{m}_{\v} p_{s,t}(x,\v)|\dif x\dif \v\leq
C(t-s)^{\frac{(\beta-n)(1+\alpha)+\gamma-m}{\alpha}},\ \ \forall s<t.\label{EV11}
\end{align}

We now use \eqref{EV11} to show the following crucial lemma, which is an analogue of Lemma \ref{Le11}.
\bl\label{Le111}
Under \eqref{CN0} and \eqref{CN00}, for any $q>-1$ and $\beta,\gamma\geq 0$ with $\beta+\gamma<\alpha$,
there is a constant $C>0$ such that for all $j\in\mN$ and $t>s\geq 0$,
\begin{align}\label{GF21}
\int^t_0\!\!\!\int_{\mR^{2d}}(t-s)^q |x|^{\beta}|\v|^{\gamma} \big|\cR^\dd_j\Gamma_{\!s,t}{p_{s,t}}(x,\v)\big|\dif x\dif\v\dif s
\leq C2^{-((1+\alpha)\beta+\gamma+(q+1)\alpha)j},
\end{align}
\begin{align}\label{GF211}
\int^t_0\!\!\!\int_{\mR^{2d}}(t-s)^q |x|^{\beta}|\v|^{\gamma} \big|\cR^x_j\Gamma_{\!s,t}{p_{s,t}}(x,\v)\big|\dif x\dif\v\dif s
\leq C2^{-(\beta+\frac{(q+1)\gamma+\alpha}{1+\alpha})j}.
\end{align}
\el
\begin{proof}
We only prove the first one. The second one is similar.
First of all, by the change of variables, we have
\begin{align*}
\sJ_{s,t}&:=\int_{\mR^{2d}}|x|^{\beta}|\v|^{\gamma} \big|\cR^\dd_j \Gamma_{\!s,t}{p_{s,t}}(x,\v)\big|\dif x\dif \v
= 2^{-((1+\alpha)(d+\beta)+d+\gamma)j}\int_{\mR^{2d}}|x|^{\beta}|\v|^{\gamma}\times\\
&\qquad\times\left|\int_{\mR^{2d}}\check\phi^\dd_1(x-\bar{x},\v-\bar{\v})p_{s,t}(2^{-(1+\alpha)j}\bar{x}+\Pi_{s,t}2^{-j}\bar\v,
2^{-j}\bar{\v})\dif \bar{x}\dif \bar{\v}\right|\dif x\dif \v.
\end{align*}
Let $\tilde U_r:=U_{(t-s)r+s}$ and $\tilde\kappa_r:=\kappa_{(t-s)r+s}$. 
By the scaling property of the heat kernel (see \cite[(2.27)]{Ch-Zh18})%\cite[(2.27)]{Ch-Zh}
,  we have
\begin{align}
p^{\kappa,U}_{s,t}(x,\v)=(t-s)^{-\frac{2d}{\alpha}-d}
p^{\tilde\kappa,\tilde U}_{0,1}((t-s)^{-\frac{1}{\alpha}-1}x,(t-s)^{-\frac{1}{\alpha}}\v).\label{NB3}
\end{align}
Hence,
$$
p^{\kappa,U}_{s,t}(2^{-(1+\alpha)j}\bar{x}+\Pi_{s,t}2^{-j}\bar\v,2^{-j}\bar{\v})=
(t-s)^{-\frac{2d}{\alpha}-d}p^{\tilde\kappa,\tilde U}_{0,1}(\hbar^{\alpha+1}\bar{x}+\hbar\theta_{s,t}\bar\v,\hbar\bar{\v}),
$$
where
$$
\hbar:=(t-s)^{-\frac{1}{\alpha}}2^{-j},\ \ \theta_{s,t}:=\Pi_{s,t}/(t-s).
$$
Since the support of $\phi^a_1$ is contained in the annulus, by Fourier's transform, 
$$
\widehat{(\Delta^{-n}_{x,\v}\check\phi^a_1)}(\xi,\eta):=(|\xi|^2+|\eta|^2)^{-n}\phi^a_1(\xi,\eta)\in\sS(\mR^{2d}),
$$ 
so that $\Delta^{-n}_{x,\v}\check\phi^a_1$ is a well-defined Schwartz function. Thus we have
\begin{align*}
\sU&:=\int_{\mR^{2d}}|x|^{\beta}|\v|^{\gamma}
\left|\int_{\mR^{2d}}\check\phi^\dd_1(x-\bar{x},\v-\bar{\v})p^{\tilde\kappa,\tilde U}_{0,1}
(\hbar^{\alpha+1}\bar x+\hbar\theta_{s,t}\bar\v,\hbar\bar \v)\dif \bar{x}\dif \bar{\v}\right|\dif x\dif \v\\
&=\int_{\mR^{2d}}|x|^{\beta}|\v|^{\gamma}\left|\int_{\mR^{2d}}\Delta^{-n}_{x,\v}\check\phi^\dd_1(x-\bar x,\v-\bar \v)\Delta^n_{x,\v}
p^{\tilde\kappa,\tilde U}_{0,1}(\hbar^{\alpha+1}\bar x+\hbar\theta_{s,t}\bar\v,\hbar\bar \v)\dif \bar{x}\dif \bar{\v}\right|\dif x\dif \v\\
&\leq\int_{\mR^{2d}}|x|^{\beta}|\v|^{\gamma}|\Delta^{-n}_{x,\v}\check\phi^\dd_1(x,\v)|\dif x\dif \v
\int_{\mR^{2d}}|\Delta^n_{x,\v}p^{\tilde\kappa,\tilde U}_{0,1}(\hbar^{\alpha+1}\bar x+\hbar\theta_{s,t}\bar\v,\hbar\bar \v)|\dif \bar{x}\dif \bar{\v}\\
&+\int_{\mR^{2d}}|\Delta^{-n}_{x,\v}\check\phi^\dd_1(x,\v)|\dif x\dif \v
\int_{\mR^{2d}}|\bar x|^{\beta}|\bar\v|^{\gamma}
|\Delta^n_{x,\v}p^{\tilde\kappa,\tilde U}_{0,1}(\hbar^{\alpha+1}\bar x+\hbar\theta_{s,t}\bar\v,\hbar\bar \v)|\dif \bar{x}\dif \bar{\v}.
\end{align*}
By the chain rule, \eqref{EV11} and cumbersome calculations, we have
$$
\int_{\mR^{2d}}|\Delta^n_{x,\v}p^{\tilde\kappa,\tilde U}_{0,1}(\hbar^{\alpha+1} x+\hbar\theta_{s,t}\v,\hbar\v)|\dif x\dif \v
\lesssim \hbar^{(\alpha+1)(n-d)-d}+\hbar^{n-(\alpha+2)d},
$$
and
$$
\int_{\mR^{2d}}|x|^{\beta}|\v|^{\gamma}
|\Delta^n_{x,\v}p^{\tilde\kappa,\tilde U}_{0,1}(\hbar^{\alpha+1}x+\hbar\theta_{s,t}\v,\hbar\v)|\dif x\dif \v
\lesssim \Big(\hbar^{(\alpha+1)(n-d)-d}+\hbar^{n-(\alpha+2)d}\Big)\hbar^{-(\alpha+1)\beta-\gamma}.
$$
Therefore,
$$
\sU\lesssim \Big(\hbar^{(\alpha+1)(n-d)-d}+\hbar^{n-(\alpha+2)d}\Big)\Big(1+\hbar^{-(\alpha+1)\beta-\gamma}\Big),
$$
and
\begin{align}\label{DA1}
\begin{split}
\sJ_{s,t}&\lesssim 2^{-((1+\alpha)\beta+\gamma)j}
\hbar^{(\alpha+2)d}\Big(\hbar^{(\alpha+1)(n-d)-d}+\hbar^{n-(\alpha+2)d}\Big)\Big(1+\hbar^{-(\alpha+1)\beta-\gamma}\Big)\\
&=\Big(\hbar^{(\alpha+1)n}+\hbar^n\Big)\Big(2^{-((1+\alpha)\beta+\gamma)j}+(t-s)^{\frac{(\alpha+1)\beta+\gamma}{\alpha}}\Big).
\end{split}
\end{align}

Without loss of generality, assume $t>2^{-\alpha j}$. We denote the left hand side of \eqref{GF21} by $\sI$, and make the following decomposition:
\begin{align*}
\sI=\left(\int^t_{t-2^{-\alpha j}}+\int^{t-2^{-\alpha j}}_0\right)(t-s)^q \sJ_{s,t}\dif s=:\sI_1+\sI_2.
\end{align*}
For $\sI_1$, using \eqref{DA1} with $n=0$, and by the change of variables, we have
\begin{align*}
\sI_1&\lesssim  \int^t_{t-2^{-\alpha j}}(t-s)^q\Big(2^{-((1+\alpha)\beta+\gamma)j}+(t-s)^{\frac{(\alpha+1)\beta+\gamma}{\alpha}}\Big)\dif s\\
&\lesssim\int^{2^{-\alpha j}}_0s^q\Big(2^{-((1+\alpha)\beta+\gamma)j}+s^{\frac{(\alpha+1)\beta+\gamma}{\alpha}}\Big)\dif s
\lesssim 2^{-((1+\alpha)\beta+\gamma+(q+1)\alpha)j}.
\end{align*}
For $\sI_2$, choosing $n$ large enough in  \eqref{DA1} so that 
$$
1+q-\tfrac{n}{\alpha}+\tfrac{(\alpha+1)\beta}{\alpha}+\tfrac{\gamma}{\alpha}<0,
$$ 
by similar calculations as above, we also have
$$
\sI_2\lesssim \int^t_{2^{-\alpha j}}s^q
\Big((s^{-\frac{1}{\alpha}}2^{-j})^{(\alpha+1)n}+(s^{-\frac{1}{\alpha}}2^{-j})^n\Big)
\Big(2^{-((1+\alpha)\beta+\gamma)j}+s^{\frac{(\alpha+1)\beta+\gamma}{\alpha}}\Big)\dif s
\lesssim 2^{-((1+\alpha)\beta+\gamma+(q+1)\alpha)j}.
$$
Combining the above calculations, we obtain the desired estimate.
\end{proof}

\section{Schauder's estimates for non-local degenerate equations}\label{SE5}

In this section we  consider the following nonlocal degenerate equation in $\mR^{2d}$:
\begin{align}\label{PDE0}
\p_t u=\sL^{(\alpha)}_{\kappa;\v} u+b\cdot\nabla u-\lambda u+f,\ \lambda\geq 0,
\end{align}
where $\sL^{(\alpha)}_{\kappa;\v}$ is defined by \eqref{KAP} and $b$ is a measurable function with the form
$$
b(t,x,\v)=(b^{(1)}(t,x,\v),b^{(2)}(t,x,\v)).
$$
Throughout this section we assume
\begin{enumerate}[{\bf (H$^{\alpha,\vartheta}_{\beta,\gamma}$)}]
\item For some $c_0\geq 1$ and $\vartheta\in(0,\alpha-1),\beta\in(0,1)$, it holds that for all $t\geq 0$ and $x,\v,w\in\mR^d$,
$$
c_0^{-1}\leq \kappa(t,x,\v,w)\leq c_0,\ \ [\kappa(t,\cdot,w)]_{\bC^\beta_\v}
+[b^{(2)}(t,\cdot)]_{\bC^{\beta}_\v}+\nor b^{(1)}(t,\cdot)\norr_{\bC^{1+\vartheta}_a}\leq c_0,
$$
where $\nor\cdot\norr_{\bC^{1+\vartheta}}$ is defined by \eqref{nor}, and for some $\gamma\in[\beta,1+\alpha)$,
$$
[\kappa(t,\cdot,w)]_{\bC^{\gamma/(1+\alpha)}_x}+[b^{(1)}(t,\cdot)]_{\bC^{(\gamma\vee(1+\vartheta))/(1+\alpha)}_x}
+[b^{(2)}(t,\cdot)]_{\bC^{\gamma/(1+\alpha)}_x}+|b(t,0)|\leq c_0,
$$
and for some closed and convex subset $\sE\subset GL_d(\mR)$, where $GL_d(\mR)$ is the set of all invertible $d\times d$-matrices,
\begin{align}\label{GA4}
\nabla_\v b^{(1)}(t,x,\v)\in \sE.
\end{align}
\end{enumerate}
%Notice that for $\gamma<\gamma'$,
%$$
%\mbox{{\bf (H$^{\beta,\vartheta}_{\gamma'}$)}}\Rightarrow \mbox{{\bf (H$^{\alpha,\vartheta}_{\beta,\gamma}$)}}.
%$$
\bd[Classical solutions]\label{Def1}
Let $\lambda\geq 0$. We call a bounded continuous function $u$ 
defined on $\mR_+\times\mR^{2d}$ a classical solution of PDE \eqref{PDE0} if for some $\eps\in(0,1)$,
$$
u\in C([0,\infty); \bC^{(\alpha\vee 1)+\eps}_\v\cap\bC^{1+\eps}_x),
$$
and for all $t\geq 0$ and $x,\v\in\mR^d$,
$$
u(t,x,\v)=\int^t_0\Big(\sL^{(\alpha)}_{\kappa;\v} u+b\cdot\nabla u-\lambda u+f\Big)(s,x,\v)\dif s.
$$
\ed

We have the following maximum principle for classical solutions.
\bt[Maximum principle]
Let $\lambda,T>0$. For any classical solution $u$ of PDE \eqref{PDE0} in the sense of Definition \ref{Def1}, it holds that
\begin{align}\label{NM4}
    \|u\|_{\mL^\infty_T}\leq (1-\e^{-\lambda T})\|f\|_{\mL^\infty_T}/\lambda.
  \end{align}
\et
\begin{proof}
Let 
\begin{align}\label{uu}
\bar u(t,x,\v):=-u(t,x,\v)\e^{\lambda t}+\int^t_0 \|f(s,\cdot,\cdot)\|_\infty \e^{\lambda s}\dif s.
\end{align}
By \eqref{PDE0}, it is easy to see that for Lebesgue almost all $t>0$,
\begin{align*}
\p_t\bar u-\sL^{(\alpha)}_{\kappa;\v} \bar u-b\cdot\nabla \bar u\geq 0.
\end{align*}
Since $\varliminf_{t \downarrow 0}\bar u(t,x,\v)=0$, by \cite[Theorem 6.1]{Ch-HXZ}, we have 
$$
\bar u(t,x,\v)\geq 0.
$$
Thus, by \eqref{uu}, we get
$$
u(t,x,\v)\leq \e^{-\lambda t}\int^t_0 \|f(s,\cdot,\cdot)\|_\infty \e^{\lambda s}\dif s\leq (1-\e^{-\lambda t})\|f\|_{\mL^\infty_t}/\lambda.
$$
By symmetry, we obtain \eqref{NM4}.
\end{proof}
The goal of this section is to prove the following Schauder's apriori estimate.
\bt\label{Th66}
Let $\alpha\in(1,2)$ and $\beta\in(0,1), \vartheta\in(0,\alpha-1)$, $\gamma\in[\beta,1+\alpha)$.
Under {\bf (H$^{\alpha,\vartheta}_{\beta,\gamma}$)}, for any $T>0$, there is a constant $C>0$ only depending on $T,c_0,\beta,\eps,d,\alpha,\sE$ such that
for any $\lambda\geq 0$ and any classical solution $u$ of \eqref{PDE0},
\begin{align}\label{Sch}
\|u\|_{\mL^\infty_T(\bC^{(\gamma+\alpha)/(1+\alpha)}_x\cap\bC^{\alpha+\beta}_\v)}
\leq C\|f\|_{\mL^\infty_T(\bC^{\gamma/(1+\alpha)}_x\cap \bC^{\beta}_\v)}.
\end{align}
\et
\br
Although our result is stated for $\alpha\in(1,2)$, it in fact also works for $\alpha=2$. In this case, under {\bf (H$^{\alpha,\beta}_{\beta,\beta}$)},
Chaudru, Honor\'e  and Menozzi \cite[Theorem 1]{Ch-Ho-Me1} has proven \eqref{Sch} for $\gamma=\beta$. When $\gamma=\beta$,
our assumption on $b^{(1)}$ is weaker since we only assume {\bf (H$^{\alpha,\vartheta}_{\beta,\beta}$)} for some $\vartheta\in(0,1)$.
\er
To prove this theorem we use the perturbation argument by freezing the coefficients along the characterization curve as usual. 
We need the following well-known fact from ODE.
\bl
Let $b:\mR_+\times\mR^{d}\to\mR^{d}$ be a time-dependent measurable vector field. Suppose that for each $t>0$, $x\mapsto b(t,x)$ is continuous and
for some $C>0$ and all $(t,x)\in\mR_+\times\mR^d$,
$$
|b(t,x)|\leq C(1+|x|).
$$
Then for each $x\in\mR^{d}$, there is a global solution $\theta_t$ to the following ODE:
$$
\dot\theta_t=b(t,\theta_t),\ \ \theta_0=x.
$$
Moreover, if we denote by $\sS_{x}:=\{\theta_\cdot: \theta_0=x\}$ the set of all solutions with starting point $x$, then for each $T>0$,
\begin{align}\label{Sur}
\cup_{x\in\mR^d}\cup_{\theta_\cdot\in\sS_{x}}\{\theta_T\}=\mR^d.
\end{align}
\el
\begin{proof}
We only show \eqref{Sur}. Fix $y\in\mR^d$ and $T>0$. Let $(\tilde\theta_t)_{t\in[0,T]}$ be the solution of ODE:
$$
\dot{\tilde\theta}_t=-b(T-t,\tilde\theta_t),\ \tilde\theta_0=y,
$$ 
and $(\bar\theta_t)_{t\geq 0}$ solve the ODE
$$
\dot{\bar\theta}_t=b(T+t,\bar\theta_t),\ \bar\theta_0=y.
$$
Define
$$
\theta_t:=\tilde\theta_{T-t} \1_{t\leq T}+\bar\theta_{t-T}\1_{t>T}.
$$
It is easy to see that $\theta_T=y$ and $\theta_\cdot\in\sS_{x}$ with $x=\tilde\theta_T$.
\end{proof}

Fix $(x_0,\v_0)\in\mR^{2d}$. Let $\theta_t$ solve the following ODE in $\mR^{2d}$:
$$
\dot\theta_t=b(t,\theta_t), \ \theta_0=(x_0,\v_0).
$$
Define
$$
\tilde u(t,x,\v):=u(t,x+\theta^{(1)}_t,\v+\theta^{ (2)}_t),\ \ \tilde f(t,x,\v):=f(t,x+\theta^{(1)}_t,\v+\theta^{ (2)}_t),
$$
$$
\kappa_0(t,w):=\kappa(t,\theta_t,w),\ \ \tilde \kappa(t, x,\v,w):=\kappa(t,x+\theta^{(1)}_t,\v+\theta^{(2)}_t,w)-\kappa(t,\theta^{(1)}_t,\theta^{(2)}_t,w),
$$
and
$$
U_t:=\nabla_\v b^{(1)}(t,\theta_t),\ \ \tilde b(t,x,\v):=b(t,x+\theta^{(1)}_t,\v+\theta^{(2)}_t)-b(t,\theta_t)-(U_t\v, 0).
$$
By \eqref{GA4}, there is a constant $c_1\geq 1$ only depending on $\sE$ such that for all $0\leq s<t$,
\begin{align}\label{UU}
|U_t|+(t-s)|\Pi_{s,t}^{-1}|\leq c_1,\mbox{ where } \Pi_{s,t}:=\int^t_s U_r\dif r.
\end{align}
It is easy to see that $\tilde u$ satisfies the following freezing equation:
\begin{align*}
\p_t \tilde u=\sL^{(\alpha)}_{\kappa_0;\v} \tilde u+U_t\v\cdot\nabla_x\tilde u-\lambda\tilde u+\sL^{(\alpha)}_{\tilde\kappa;\v} \tilde u+\tilde b\cdot\nabla\tilde u+\tilde f,
\end{align*}
where
$$
\sL^{(\alpha)}_{\kappa_0;\v} u(x,\v):=\int_{\mR^d}\delta^{(2)}_{w} u(x,\v)\kappa_0(t,w)\frac{\dif w}{|w|^{d+\alpha}},\ \ \delta^{(2)}_w:=\delta_{w;2}^*\delta_{w;2}.
$$
Below, without loss of generality, we drop the tilde over $u,f,\kappa,b$ and assume $x_0=\v_0=0$ and
\begin{align}\label{EM1}
|\kappa(t,x,\v,w)|\leq [\kappa(t,\cdot,w)]_{\bC^{\gamma/(1+\alpha)}_x}|x|^{\frac{\gamma}{1+\alpha}}+[\kappa(t,\cdot,w)]_{\bC^{\beta}_x}|\v|^{\beta},
\end{align}
and
% under {\bf (H$^{\beta}_b$)},
\begin{align}\label{EM2}
&\quad|b^{(1)}(t,x,\v)|\leq \nor b^{(1)}(t,\cdot)\norr_{\bC^{1+\vartheta}_\dd}\big(|x|^{\frac{1+\vartheta}{1+\alpha}}+|\v|^{1+\vartheta}\big),\\
&|b^{(2)}(t,x,\v)|\leq [b^{(2)}(t,\cdot)]_{\bC^{\gamma/(1+\alpha)}_x}|x|^{\frac{\gamma}{1+\alpha}}+[b^{(2)}(t,\cdot)]_{\bC^{\beta}_\v}|\v|^{\beta}.\label{EM22}
\end{align}
Let $p^{\kappa_0}_{s,t}(x,\v)$ be the heat kernel of $\sL^{(\alpha)}_{\kappa_0;\v}+U_t\v\cdot\nabla_x$. Define for $\lambda\geq 0$,
\begin{align}\label{EM98}
P^\lambda_{s,t}f(x,\v):=\Big(\Gamma_{\!s,t}p^\lambda_{s,t}*\Gamma_{\!s,t}f\Big)(x,\v),\quad 
p^\lambda_{s,t}(x,\v):=\e^{\lambda(s-t)}p^{\kappa_0}_{s,t}(x,\v),
\end{align}
where $\Gamma_{\!s,t}$ is defined by \eqref{DA4}.
By Duhamel's formula we have (see \eqref{DA3})
\begin{align}\label{EM8}
\begin{split}
u(t,x,\v)&=\int^t_0P^\lambda_{s,t}\sL^{(\alpha)}_{\kappa;\v} u(s,x,\v)\dif s
+\int^t_0P^\lambda_{s,t}(b\cdot\nabla u)(s,x,\v)\dif s
+\int^t_0P^\lambda_{s,t} f(s,x,\v)\dif s.
\end{split}
\end{align}

We prepare the following lemmas.
\bl\label{Le522}
Under {\bf (H$^{\alpha,\vartheta}_{\beta,\gamma}$)}, for any $\eps\in(0,1)$, 
there is a constant $C>0$ such that for all $j\in\mN$ and $\lambda\geq 0$, $t\in[0,T]$,
\begin{align}\label{DA40}
\int^t_0|\cR^a_jP^\lambda_{s,t}\sL^{(\alpha)}_{\kappa;\v} u|(s,0,0)\dif s&\leq C2^{-(\alpha+\beta)j}\|u\|_{\mL^\infty_T(\bC^{\alpha+\eps}_{\dd})},\\
\int^t_0|\cR^x_jP^\lambda_{s,t}\sL^{(\alpha)}_{\kappa;\v} u|(s,0,0)\dif s&\leq C2^{-\frac{\gamma+\alpha}{1+\alpha}j}
\|u\|_{\mL^\infty_T(\bC^{(\gamma-\beta)/(1+\alpha)}_x\bC^{\alpha+\eps}_\v)}.\label{DA425}
\end{align}
\el
\begin{proof}
(i) First of all, by \eqref{EM1}, we have for $u\in\bC^{\alpha+\eps}_\v$,
\begin{align*}
\left|\sL^{(\alpha)}_{\kappa;\v} u(x,\v)\right|
&=\left|\int_{\mR^d}\delta^{(2)}_{w} u(x,\v)\cdot \kappa(t,x,\v,\w)\frac{\dif w}{|w|^{d+\alpha}}\right|
\lesssim\Big(|x|^{\frac{\gamma}{1+\alpha}}+|\v|^{\beta}\Big)\int_{\mR^d}\big|\delta^{(2)}_{w} u(x,\v)\big|\frac{\dif w}{|w|^{d+\alpha}}\\
&\lesssim\Big(|x|^{\frac{\gamma}{1+\alpha}}+|\v|^{\beta}\Big)\|u\|_{\bC^{\alpha+\eps}_\v}\int_{\mR^d}(1\wedge|w|^{\alpha+\eps})\frac{\dif w}{|w|^{d+\alpha}}
\lesssim\Big(|x|^{\frac{\gamma}{1+\alpha}}+|\v|^{\beta}\Big)\|u\|_{\bC^{\alpha+\eps}_\v}.
\end{align*}
Thus by definition \eqref{EM98}, we have
\begin{align*}
\int^t_0|\cR^a_jP^\lambda_{s,t}\sL^{(\alpha)}_{\kappa;\v} u|(s,0,0)\dif s&=\int^t_0\left|\int_{\mR^{2d}}\cR^\dd_j\Gamma_{\!s,t} p^\lambda_{s,t}(x,\v)
\cdot\Big(\Gamma_{\!s,t}\sL^{(\alpha)}_{\kappa;\v}  u\Big)(s,x,\v)\dif x\dif \v\right|\dif s\\
&\lesssim \|u\|_{\mL^\infty_T(\bC^{\alpha+\eps}_\v)}\int^t_0\!\!\!\int_{\mR^{2d}}|\cR^\dd_j\Gamma_{\!s,t} p^\lambda_{s,t}(x,\v)|
\Big(|x+\Pi_{s,t}\v|^{\frac{\gamma}{1+\alpha}}+|\v|^{\beta}\Big)\dif x\dif \v\dif s,
\end{align*}
which in turn gives \eqref{DA40}  by direct application of  \eqref{GF21} and $\gamma\geq\beta$.

(ii) Notice that by \eqref{KJ2} and $\cR^x_j\Gamma_{\!s,t}=\Gamma_{\!s,t}\cR^x_j$,
\begin{align*}
&\int^t_0|\cR^x_jP^\lambda_{s,t}\sL^{(\alpha)}_{\kappa;\v} u|(s,0,0)\dif s
=\int^t_0\left|\int_{\mR^{2d}}\cR^x_j\Gamma_{\!s,t} p^\lambda_{s,t}(x,\v)
\cdot\Big(\Gamma_{\!s,t}\sL^{(\alpha)}_{\kappa;\v}  u\Big)(s,x,\v)\dif x\dif \v\right|\dif s\\
&\qquad=\int^t_0\left|\int_{\mR^{2d}}\widetilde\cR^x_j\Gamma_{\!s,t} p^\lambda_{s,t}(x,\v)
\cdot\Big(\Gamma_{\!s,t}\cR^x_j\sL^{(\alpha)}_{\kappa;\v}  u\Big)(s,x,\v)\dif x\dif \v\right|\dif s=:\sI_1+\sI_2,
\end{align*}
where
\begin{align*}
\sI_1&:=\int^t_0\left|\int_{\mR^{2d}}\widetilde\cR^x_j\Gamma_{\!s,t} p^\lambda_{s,t}(x,\v)
\cdot\Big(\Gamma_{\!s,t}\sL^{(\alpha)}_{\kappa;\v}\cR^x_j  u\Big)(s,x,\v)\dif x\dif \v\right|\dif s,\\
\sI_2&:=\int^t_0\left|\int_{\mR^{2d}}\widetilde\cR^x_j\Gamma_{\!s,t} p^\lambda_{s,t}(x,\v)
\cdot\Big(\Gamma_{\!s,t}[\cR^x_j,\sL^{(\alpha)}_{\kappa;\v}]  u\Big)(s,x,\v)\dif x\dif \v\right|\dif s.
\end{align*}
For $\sI_1$, by the assumptions, we have
\begin{align*}
\left|\sL^{(\alpha)}_{\kappa;\v}\cR^x_j u(x,\v)\right|
&=\left|\int_{\mR^d}\delta^{(2)}_{w}\cR^x_j u(x,\v)\cdot \kappa(t,x,\v,\w)\frac{\dif w}{|w|^{d+\alpha}}\right|\\
&\lesssim\Big(|x|^{\frac{\gamma}{1+\alpha}}+|\v|^{\beta}\Big)\int_{\mR^d}\big|\delta^{(2)}_{w} \cR^x_ju(x,\v)\big|\frac{\dif w}{|w|^{d+\alpha}}\\
&\lesssim\Big(|x|^{\frac{\gamma}{1+\alpha}}+|\v|^{\beta}\Big)\|\cR^x_ju\|_{\bC^{\alpha+\eps}_\v},
\end{align*}
and thus, by \eqref{GF211} and $\gamma\geq\beta$,
\begin{align*}
\sI_1&\lesssim\|\cR^x_ju\|_{\mL^\infty_T(\bC^{\alpha+\eps}_\v)}\int^t_0\int_{\mR^{2d}}|\widetilde\cR^x_j\Gamma_{\!s,t} p^\lambda_{s,t}|(x,\v)
\Big(|x+\Pi_{s,t}\v|^{\frac{\gamma}{1+\alpha}}+|\v|^{\beta}\Big)\dif x\dif \v\dif s\\
&\lesssim\|\cR^x_ju\|_{\mL^\infty_T(\bC^{\alpha+\eps}_\v)} \Big(2^{-\frac{\gamma+\alpha}{1+\alpha}j}+2^{-\frac{\alpha+\beta}{1+\alpha}j}\Big)
\lesssim 2^{-\frac{\gamma+\alpha}{1+\alpha}j}\|u\|_{\mL^\infty_T(\bC^{(\gamma-\beta)/(1+\alpha)}_x\bC^{\alpha+\eps}_\v)}.
\end{align*}
For $\sI_2$, by definition, we have
\begin{align*}
\|[\cR^x_j,\sL^{(\alpha)}_{\kappa;\v}]  u\|_\infty\lesssim
\int_{\mR^{2d}}\|\delta^{(2)}_{w} u\|_\infty
\frac{\dif w}{|w|^{d+\alpha}}\int_{\mR^d}|\bar x|^{\frac{\gamma}{1+\alpha}}\check\phi_j( \bar x)\dif \bar x
\lesssim \|u\|_{\bC^{\alpha+\eps}_\v}2^{-\frac{\gamma}{1+\alpha}j}.
\end{align*}
Hence,  by \eqref{GF211} again,
\begin{align*}
\sI_2\lesssim 
\|u\|_{\mL^\infty_T(\bC^{\alpha+\eps}_\v)}2^{-\frac{\gamma}{1+\alpha}j}\int^t_0\!\!\int_{\mR^{2d}}|\widetilde\cR^x_j\Gamma_{\!s,t} p^\lambda_{s,t}|(x,\v)\dif x\dif \v\dif s
\lesssim \|u\|_{\mL^\infty_T(\bC^{\alpha+\eps}_\v)}2^{-\frac{\gamma+\alpha}{1+\alpha}j}.
\end{align*}
Combining the above calculations, we obtain \eqref{DA425}.
\end{proof}

To treat the other terms in \eqref{EM8}, we need the following lemma.
\bl\label{Le53}
Let $c_1\geq 1$ be the same as in \eqref{UU}. For $t\geq 0$ and $j\in\mN$, define
\begin{align*}
\Theta^{t}_j:=\Big\{\ell\in\mN_0: 2^{\ell}\leq 2^4 c_1(2^j+t2^{(1+\alpha)j}),\ 
2^{j}\leq 2^4 c_1(2^\ell+t2^{(1+\alpha)\ell})\Big\}.
\end{align*}
\begin{enumerate}[(i)]
\item Let $0\leq s<t$ and $j\in\mN$. For any $\ell\notin\Theta^{t-s}_j$, it holds that
\begin{align}\label{EM3}
\<\cR^\dd_jf,\Gamma_{\!s,t}\cR^\dd_\ell g\>=\int_{\mR^{2d}}\cR^\dd_jf(x,\v)\cdot\Gamma_{\!s,t}
\cR^\dd_\ell g(x,\v)\dif x\dif\v=0.
\end{align}
\item For any $\beta>0$, there is a constant $C=C(c_1,\beta)>0$ such that for all $j\in\mN$ and $t\geq 0$,
\begin{align}\label{DA2}
\sum_{\ell\in\Theta^t_j}2^{-\beta\ell}\leq C\Big(2^{-j}+t2^{(\alpha-1)j}\Big)^{\beta},\ \ \sum_{\ell\in\Theta^t_j}2^{\beta\ell}\leq 
C\Big(2^j+t2^{(1+\alpha)j}\Big)^\beta.
\end{align}
\item For any $T>0$, there exists a $j_0=j_0(c_1,\alpha,T)\in\mN$ such that for all $j\geq j_0$ and $t\in[0,T]$,
$$
\inf\Big\{\ell: \ell\in \Theta^{t}_j\Big\}\geq 5.
$$
\end{enumerate}
\el
\begin{proof}
(i) To prove \eqref{EM3}, by Fourier's transform we have
\begin{align*}
\<\cR^\dd_jf,\Gamma_{\!s,t}\cR^\dd_\ell g\>
=\int_{\mR^{2d}}\phi^\dd_j(\xi,\eta)\hat{f}(\xi,\eta)\phi^\dd_\ell(\xi,\eta-\Pi_{s,t}\xi)\hat{g}(\xi,\eta-\Pi_{s,t}\xi)\dif \xi\dif \eta.
\end{align*}
Notice that
$$
\mbox{supp}\phi^\dd_j\subset\Big\{(\xi,\eta): 2^{j-1}\leq |\xi|^{1/(1+\alpha)}+|\eta|\leq 2^{j+1}\Big\}=:\cI_j.
$$
Assuming $\<\cR^\dd_j\Gamma_{\!s,t}f,\Gamma_{\!s,t}\cR^\dd_\ell g\>\not=0$ for $j,\ell\in\mN$, we must have
$$
(\xi,\eta)\in\cI_j\ \mbox{ and }\ (\xi,\eta-\Pi_{s,t}\xi)\in\cI_\ell,
$$
which implies that
$$
|\xi|\leq 2^{(1+\alpha)(j+1)},\ \ |\eta|\leq 2^{j+1},
$$
and
$$
2^{\ell-1}\leq|\xi|^{1/(1+\alpha)}+|\eta-\Pi_{s,t}\xi|
\leq 2\cdot 2^{j+1}+c_1(t-s) 2^{(1+\alpha)(j+1)}\leq 2^3c_1(2^j+(t-s)2^{(1+\alpha)j}).
$$
By symmetry we also have
$$
2^{j-1}\leq 2^3 c_1(2^\ell+(t-s)2^{(1+\alpha)\ell}).
$$
If $\<\cR^\dd_j\Gamma_{\!s,t}f,\Gamma_{\!s,t}\cR^\dd_0 g\>\not=0$ for $j\in\mN$, we still have 
$$
2^{j-1}\leq |\xi|^{1/(1+\alpha)}+|\eta-\Pi_{s,t}\xi|+|\Pi_{s,t}\xi|\leq 2+c_1(t-s) 2^{1+\alpha}.
$$
Combining the above calculations, one sees that for $\ell\notin\Theta^{t-s}_j$, \eqref{EM3} holds.

(ii) We only prove the first estimate in \eqref{DA2}. If $\ell>j$, then $2^{-\ell}\leq 2^{-j}$. If $\ell\leq j$, then by the definition of $\Theta^t_j$, 
$$
2^{-\ell}\leq 2^4 c_12^{-j}(1+(t-s)2^{\alpha\ell})\leq 2^4 c_1(2^{-j}+(t-s)2^{(\alpha-1)j})=:D,
$$
which implies $\ell\geq -\ln D/\ln 2$. Thus, we have
$$
\sum_{\ell\in\Theta^t_j}2^{-\beta\ell}\leq\sum_{\ell\geq -\ln D/\ln 2}2^{-\beta\ell}\leq (2 D)^\beta/(1-2^{-\beta})\lesssim (2^{-j}+(t-s)2^{(\alpha-1)j})^\beta.
$$

(iii) By definition of $\Theta^{t}_j$, it suffices to take $j_0>\ln(2^4 c_1(2^5+T2^{(1+\alpha)5}))/\ln 2$.
\end{proof}

\bl\label{Le54}
Let $T>0$ and $j_0$ be as in (iii) of Lemma \ref{Le53}.
Under {\bf (H$^{\alpha,\vartheta}_{\beta,\gamma}$)},
there is a constant $C>0$ such that for all $j\geq j_0$, $\lambda\geq 0$ and $t\in[0,T]$,
\begin{align}\label{DA5}
\int^t_0|\cR^a_j P^\lambda_{s,t}(b\cdot\nabla u)|(s,0,0)\dif s&\leq C2^{-(\alpha+\beta)j}\|u\|_{\mL^\infty_T(\bC^{\alpha+\beta-\vartheta}_{\dd})},\\
\int^t_0|\cR^x_jP^\lambda_{s,t}(b^{(1)}\cdot\nabla_x u)|(s,0,0)\dif s&\leq C2^{-\frac{\gamma+\alpha}{1+\alpha}j}
\|u\|_{\mL^\infty_T(\bC^{(\gamma+\alpha-\vartheta)/(1+\alpha)}_x)},\label{DA405}\\
\int^t_0|\cR^x_jP^\lambda_{s,t}(b^{(2)}\cdot\nabla_\v u)|(s,0,0)\dif s&\leq C2^{-\frac{\gamma+\alpha}{1+\alpha}j}
\|\nabla_\v u\|_{\mL^\infty_T(\bC^{(\gamma-\beta)/(1+\alpha)}_x)}.\label{DA4035}
\end{align}
\el
\begin{proof}
(i) Let $\Theta=\Theta^{t-s}_j$ be as in Lemma \ref{Le53}. By definition \eqref{EM98} and \eqref{EM3}, we have
\begin{align}
&\int^t_0|\cR^a_j P^\lambda_{s,t}(b\cdot\nabla u)|(s,0,0)\dif s
=\int^t_0\left|\int_{\mR^{2d}}\cR^\dd_j\Gamma_{\!s,t} p^\lambda_{s,t}(x,\v)
\cdot\Gamma_{\!s,t}(b\cdot\nabla u)(s,x,\v)\dif x\dif \v\right|\dif s\no\\
&\qquad\qquad=\int^t_0\left|\sum_{\ell\in\Theta}\int_{\mR^{2d}}\cR^\dd_j\Gamma_{\!s,t} p^\lambda_{s,t}(x,\v)
\cdot\Gamma_{\!s,t}\cR^\dd_\ell(b\cdot\nabla u)(s,x,\v)\dif x\dif \v\right|\dif s.\label{GZ1}
\end{align}
Below we drop the time variable and write
\begin{align*}
\cR^\dd_\ell(b\cdot\nabla u)(x,\v)=(b\cdot\nabla\cR^\dd_\ell u)(x,\v)+
[\cR^\dd_\ell, b\cdot\nabla] u(x,\v)=:\sI_1+\sI_2.
\end{align*}
For $\sI_1$,  by \eqref{EM2} and \eqref{GA6}, we have
\begin{align*}
|(b^{(1)}\cdot\nabla_x\cR^\dd_\ell u)(x,\v)|&\lesssim 
\Big(|x|^{\frac{1+\vartheta}{1+\alpha}}+|\v|^{1+\vartheta}\Big) 2^{(1+\alpha)\ell}\|\cR^\dd_\ell u\|_\infty
\lesssim\Big(|x|^{\frac{1+\vartheta}{1+\alpha}}+|\v|^{1+\vartheta}\Big) 2^{(1-\beta+\vartheta)\ell}\|u\|_{\bB^{\alpha+\beta-\vartheta}_{\dd,\infty}},
\end{align*}
and by \eqref{EM22} and \eqref{GA6},
\begin{align*}
|(b^{(2)}\cdot\nabla_\v\cR^\dd_\ell u)(x,\v)|&\lesssim\Big(|x|^{\frac{\gamma}{1+\alpha}}+|\v|^{\beta}\Big)2^\ell\|\cR^\dd_\ell u\|_\infty
\lesssim\Big(|x|^{\frac{\gamma}{1+\alpha}}+|\v|^{\beta}\Big)2^{(1-\alpha-\beta+\vartheta)\ell}\|u\|_{\bB^{\alpha+\beta-\vartheta}_{\dd,\infty}}.
\end{align*}
Hence, by \eqref{LJ1},
\begin{align}\label{GZ2}
|\sI_1|\lesssim 
\|u\|_{\bC^{\alpha+\beta-\vartheta}_{\dd}}\Big(\big(|x|^{\frac{1+\vartheta}{1+\alpha}}+|\v|^{1+\vartheta}\big) 2^{(1-\beta+\vartheta)\ell}
+\big(|x|^{\frac{\gamma}{1+\alpha}}+|\v|^{\beta}\big)2^{(1-\alpha-\beta+\vartheta)\ell}\Big).
\end{align}
For $\sI_2$, due to $j\geq j_0$ and by (iii) of Lemma \ref{Le53}, we have $\ell\geq 5$ for $\ell\in\Theta^{t-s}_j$.
Thus we can use \eqref{HQ2} with $(\vartheta,\beta-1-\vartheta)$ in place of $(\beta,\gamma)$ to derive that
\begin{align*}
\big|\big[\cR^\dd_\ell,b^{(1)}\cdot\nabla_x\big]  u(x,\v)\big|
&\lesssim2^{-\ell(\beta-\vartheta)}\Big(2^{-\ell\vartheta}+|x|^{\frac{\vartheta}{1+\alpha}}+|\v|^{\vartheta}\Big)
\|\nabla_x u\|_{\bC^{\beta-1-\vartheta}_{\dd}}\\
&\lesssim\Big(2^{-\ell \beta}+2^{-\ell(\beta-\vartheta)}(|x|^{\frac{\vartheta}{1+\alpha}}+|\v|^{\vartheta})\Big)
\|u\|_{\bC^{\alpha+\beta-\vartheta}_{\dd}}.
\end{align*}
Moreover, by $\gamma\geq\beta$, $\alpha-\vartheta>1$ and the definition, we also have
\begin{align*}
\big|\big[\cR^\dd_\ell,b^{(2)}\cdot\nabla_\v\big]  u(x,\v)\big|
&=\left|\int_{\mR^{2d}}\check\phi^\dd_\ell(x-\bar x,\v-\bar \v)
\Big(b^{(2)}(\bar x,\bar \v)-b^{(2)}(x,\v)\Big)\nabla_\v u(\bar x,\bar\v)\dif \bar x\dif\bar\v\right|\\
&\lesssim\|\nabla_\v u\|_\infty
\int_{\mR^{2d}}|\check\phi^\dd_\ell(\bar x,\bar \v)|\cdot\Big(|\bar x|^{\frac{\gamma}{1+\alpha}}+|\bar\v|^\beta)\dif \bar x\dif\bar\v
\lesssim 2^{-\ell\beta}\|u\|_{\bC^{\alpha+\beta-\vartheta}_{\dd}}.
\end{align*}
Therefore,
\begin{align}\label{GZ3}
|\sI_2|\lesssim \|u\|_{\bC^{\alpha+\beta-\vartheta}_{\dd}}\Big(2^{-\ell \beta}+2^{-\ell(\beta-\vartheta)}(|x|^{\frac{\vartheta}{1+\alpha}}+|\v|^{\vartheta})\Big).
\end{align}
Combining \eqref{GZ1}-\eqref{GZ3}, and by \eqref{DA2} we get
\begin{align*}
&\int^t_0|\cR^a_j P^\lambda_{s,t}(b\cdot\nabla u)|(s,0,0)\dif s
\lesssim \|u\|_{\mL^\infty_T(\bC^{\alpha+\beta-\vartheta}_{\dd})}
\int^t_0\!\!\!\int_{\mR^{2d}}\left|\cR^\dd_j\Gamma_{\!s,t} p^\lambda_{s,t}(x,\v)\right|\\
&\qquad\qquad\times\Big\{\big(|x+\Pi_{s,t}\v|^{\frac{1+\vartheta}{1+\alpha}}+|\v|^{1+\vartheta}\big)
\big(2^{j}+(t-s)2^{(1+\alpha)j}\big)^{1-\beta+\vartheta}\\
&\qquad\qquad\quad+\big(|x+\Pi_{s,t}\v|^{\frac{\gamma}{1+\alpha}}+|\v|^{\beta}\big)
\big(2^{-j}+(t-s)2^{(\alpha-1)j}\big)^{\alpha+\beta-\vartheta-1}\\
&\qquad\qquad\quad+\big(2^{-j}+(t-s)2^{(\alpha-1)j}\big)^{\beta-\vartheta}
\big(|x+\Pi_{s,t}\v|^{\frac{\vartheta}{1+\alpha}}+|\v|^{\vartheta}\big)\\
&\qquad\qquad\quad+\big(2^{-j}+(t-s)2^{(\alpha-1)j}\big)^{\beta}\Big\}\dif x\dif\v\dif s,
\end{align*}
which in turn yields \eqref{DA5} by using \eqref{UU} and \eqref{GF21} item by item, as well as $\gamma\geq\beta$.

(ii) As above by definition \eqref{EM98} and \eqref{KJ2}, we have
\begin{align*}
&\int^t_0|\cR^x_jP^\lambda_{s,t}(b^{(1)}\cdot\nabla_x u)|(s,0,0)\dif s
=\int^t_0\left|\int_{\mR^{2d}}\cR^x_j\Gamma_{\!s,t} p^\lambda_{s,t}(x,\v)
\cdot\Gamma_{\!s,t}(b^{(1)}\cdot\nabla_x u)(s,x,\v)\dif x\dif \v\right|\dif s\\
&\qquad=\int^t_0\left|\int_{\mR^{2d}}\widetilde\cR^x_j\Gamma_{\!s,t} p^\lambda_{s,t}(x,\v)
\cdot\Gamma_{\!s,t}\cR^x_j(b^{(1)}\cdot\nabla_x u)(s,x,\v)\dif x\dif \v\right|\dif s=\sI_1+\sI_2,
\end{align*}
where we have used $\cR^x_j\Gamma_{\!s,t}=\Gamma_{\!s,t}\cR^x_j$, and
\begin{align*}
\sI_1&:=\int^t_0\left|\int_{\mR^{2d}}\widetilde\cR^x_j\Gamma_{\!s,t} p^\lambda_{s,t}(x,\v)
\cdot\Gamma_{\!s,t}(b^{(1)}\cdot\nabla_x \cR^x_j u)(s,x,\v)\dif x\dif \v\right|\dif s,\\
\sI_2&:=\int^t_0\left|\int_{\mR^{2d}}\widetilde\cR^x_j\Gamma_{\!s,t} p^\lambda_{s,t}(x,\v)
\cdot\Gamma_{\!s,t}[\cR^x_j,b^{(1)}\cdot\nabla_x]  u(s,x,\v)\dif x\dif \v\right|\dif s.
\end{align*}
For $\sI_1$, noticing that by \eqref{EM2} and the definition of $\cR^x_j$,
\begin{align*}
|b^{(1)}\cdot\nabla_x \cR^x_j u|(x,\v)\lesssim \Big(|x|^{\frac{1+\vartheta}{1+\alpha}}+|\v|^{1+\vartheta}\Big) 2^{j}\|\cR^x_j u\|_\infty,
\end{align*}
we have
\begin{align*}
\sI_1&\lesssim 2^{j}\|\cR^x_j u\|_{\mL^\infty_T}
\int^t_0\!\!\!\int_{\mR^{2d}}|\widetilde\cR^x_j\Gamma_{\!s,t} p^\lambda_{s,t}|(x,\v)\Big(|x+\Pi_{s,t}\v|^{\frac{1+\vartheta}{1+\alpha}}
+|\v|^{1+\vartheta}\Big) \dif x\dif \v\dif s\\
&\stackrel{\eqref{GF211}}{\lesssim}\|\cR^x_j u\|_{\mL^\infty_T} 2^{-\frac{\vartheta}{1+\alpha}j}
\lesssim\|u\|_{\mL^\infty_T(\bB^{(\gamma+\alpha-\vartheta)/(1+\alpha)}_{x,\infty})} 
2^{-\frac{\gamma+\alpha}{1+\alpha}j}.
\end{align*}
For $\sI_2$, by \eqref{GS1} with $(\frac{(1+\vartheta)\vee\gamma}{1+\alpha},\frac{-((1+\vartheta)\vee\gamma-\gamma)}{1+\alpha})$ in place of $(\beta,\gamma)$
and \eqref{GF211}, we have
\begin{align*}
\sI_2&\lesssim 
\int^t_0\int_{\mR^{2d}}|\widetilde\cR^x_j\Gamma_{\!s,t} p^\lambda_{s,t}|(x,\v)\,
\|[\cR^x_j,b^{(1)}\cdot\nabla_x]  u\|_\infty\dif x\dif \v\dif s\\
&\lesssim 2^{-\frac{\gamma+\alpha}{1+\alpha}j}\|\nabla_x u\|_{\mL^\infty_T(\bC^{-((1+\vartheta)\vee\gamma-\gamma)/(1+\alpha)}_{x})}
\lesssim 2^{-\frac{\gamma+\alpha}{1+\alpha}j}\|u\|_{\mL^\infty_T(\bC^{(\gamma+\alpha-\vartheta)/(1+\alpha)}_{x})}.
\end{align*}
Combining the above two estimates and Theorem \ref{Th22}, we obtain \eqref{DA405}.

(iii) As above by definition we have
\begin{align*}
&\int^t_0|\cR^x_jP^\lambda_{s,t}(b^{(2)}\cdot\nabla_\v u)|(s,0,0)\dif s
=\int^t_0\left|\int_{\mR^{2d}}\cR^x_j\Gamma_{\!s,t} p^\lambda_{s,t}(x,\v)
\cdot\Gamma_{\!s,t}(b^{(2)}\cdot\nabla_\v u)(s,x,\v)\dif x\dif \v\right|\dif s\\
&\qquad=\int^t_0\left|\int_{\mR^{2d}}\widetilde\cR^x_j\Gamma_{\!s,t} p^\lambda_{s,t}(x,\v)
\cdot\Gamma_{\!s,t}\cR^x_j(b^{(2)}\cdot\nabla_\v u)(s,x,\v)\dif x\dif \v\right|\dif s=\sI_1+\sI_2,
\end{align*}
where
\begin{align*}
\sI_1&:=\int^t_0\left|\int_{\mR^{2d}}\widetilde\cR^x_j\Gamma_{\!s,t} p^\lambda_{s,t}(x,\v)
\cdot\Gamma_{\!s,t}(b^{(2)}\cdot\nabla_\v \cR^x_j u)(s,x,\v)\dif x\dif \v\right|\dif s,\\
\sI_2&:=\int^t_0\left|\int_{\mR^{2d}}\widetilde\cR^x_j\Gamma_{\!s,t} p^\lambda_{s,t}(x,\v)
\cdot\Gamma_{\!s,t}[\cR^x_j,b^{(2)}\cdot\nabla_\v]  u(s,x,\v)\dif x\dif \v\right|\dif s.
\end{align*}
For $\sI_1$, noticing that by \eqref{EM22},
\begin{align*}
|b^{(2)}\cdot\nabla_\v \cR^x_j u|(x,\v)\lesssim \Big(|x|^{\frac{\gamma}{1+\alpha}}+|\v|^{\beta}\Big) \|\nabla_\v\cR^x_j u\|_\infty,
\end{align*}
we have
\begin{align*}
\sI_1&\lesssim \|\nabla_\v\cR^x_j u\|_{\mL^\infty_T}
\int^t_0\!\!\!\int_{\mR^{2d}}|\widetilde\cR^x_j\Gamma_{\!s,t} p^\lambda_{s,t}|(x,\v)\Big(|x+\Pi_{s,t}\v|^{\frac{\gamma}{1+\alpha}}+|\v|^{\beta}\Big) \dif x\dif \v\dif s\\
&\stackrel{\eqref{GF211}}{\lesssim} \|\nabla_\v\cR^x_j u\|_{\mL^\infty_T}\Big(2^{-\frac{\gamma+\alpha}{1+\alpha}j}+2^{-\frac{\alpha+\beta}{1+\alpha}j}\Big)
\lesssim 2^{-\frac{\gamma+\alpha}{1+\alpha}j}\|\nabla_\v u\|_{\mL^\infty_T(\bB^{(\gamma-\beta)/(1+\alpha)}_{x,\infty})}.
\end{align*}
For $\sI_2$, by the commutator estimate \eqref{GS1}, we have
\begin{align*}
\sI_2&\lesssim 
\int^t_0\int_{\mR^{2d}}|\widetilde\cR^x_j\Gamma_{\!s,t} p^\lambda_{s,t}|(x,\v)
\|[\cR^x_j,b^{(2)}\cdot\nabla_\v]  u\|_\infty\dif x\dif \v\dif s\lesssim 2^{-\frac{\gamma+\alpha}{1+\alpha}j}\|\nabla_\v u\|_{\mL^\infty_T}.
\end{align*}
Combining the above calculations, we obtain \eqref{DA4035}.
\end{proof}

\bl\label{Le58}
For any $\beta\in(0,1)$,
there is a constant $C>0$ such that for all $j\geq 5$ and $\lambda\geq 0$,  $t\in[0,T]$,
\begin{align}
\int^t_0|\cR^a_jP^\lambda_{s,t}f|(s,0,0)\dif s&\leq C2^{-(\alpha+\beta)j}\|f\|_{\mL^\infty_T(\bC^{\beta}_{\dd})},\label{DA44}\\
\int^t_0|\cR^x_jP^\lambda_{s,t}f|(s,0,0)\dif s&\leq C2^{-\frac{\gamma+\alpha}{1+\alpha}j}\|f\|_{\mL^\infty_T(\bC^{\gamma/(1+\alpha)}_x)}.\label{DA494}
\end{align}
\el
\begin{proof}
We only prove the first one. The second one is similar and easier by \eqref{GF211}.
Let $\Theta=\Theta^{t-s}_j$  be as in Lemma \ref{Le53}. By definition \eqref{EM98} and Lemma \ref{Le53}, we have
\begin{align*}
&\int^t_0|\cR^a_j P^\lambda_{s,t}f|(s,0,0)\dif s
=\int^t_0\left|\int_{\mR^{2d}}\cR^\dd_j\Gamma_{\!s,t} p^\lambda_{s,t}(x,\v)
\cdot\Gamma_{\!s,t}f(s,x,\v)\dif x\dif \v\right|\dif s\\
&\qquad=\int^t_0\left|\sum_{\ell\in\Theta}\int_{\mR^{2d}}\cR^\dd_j\Gamma_{\!s,t} p^\lambda_{s,t}(x,\v)
\cdot\Gamma_{\!s,t}\cR^\dd_\ell f(s,x,\v)\dif x\dif \v\right|\dif s\\
&\qquad\leq\int^t_0\sum_{\ell\in\Theta}\|\cR^\dd_\ell f(s)\|_\infty\left(\int_{\mR^{2d}}|\cR^\dd_j\Gamma_{\!s,t} p^\lambda_{s,t}(x,\v)|\dif x\dif \v\right)\dif s\\
&\qquad\leq \|f\|_{\mL^\infty_T(\bB^{\beta}_{\dd,\infty})}\int^t_0\sum_{\ell\in\Theta}2^{-\ell\beta}
\left(\int_{\mR^{2d}}|\cR^\dd_j\Gamma_{\!s,t} p^\lambda_{s,t}(x,\v)|\dif x\dif \v\right)\dif s\\
&\qquad\lesssim \|f\|_{\mL^\infty_T(\bC^{\beta}_{\dd})}\int^t_0\Big(2^{-j}+(t-s)2^{(\alpha-1)j}\Big)^{\beta}
\left(\int_{\mR^{2d}}|\cR^\dd_j\Gamma_{\!s,t} p^\lambda_{s,t}(x,\v)|\dif x\dif \v\right)\dif s,
\end{align*}
which gives \eqref{DA44}  by  application of  \eqref{GF21}.
\end{proof}

Now we are in a position to give
\begin{proof}[Proof of Theorem \ref{Th66}]

(i) Fix $\eps\in(0,1)$. We first show the following estimates:
\begin{align}\label{Sch1}
\|u\|_{\mL^\infty_T(\bC^{\alpha+\beta}_{\dd})}
\leq C\|f\|_{\mL^\infty_T(\bC^\beta_{\dd})},
\end{align}
and
\begin{align}\label{FW1}
\|u\|_{\mL^\infty_T(\bC^{(\gamma+\alpha)/(1+\alpha)}_x)}
\leq C\Big(\|f\|_{\mL^\infty_T(\bC^{\gamma/(1+\alpha)}_x)}+\|u\|_{\mL^\infty_T(\bC^{(\gamma-\beta)/(1+\alpha)}_x\bC^{\alpha+\eps}_\v)}\Big).
\end{align}
By Lemmas \ref{Le522}, \ref{Le54} and \ref{Le58}, we have
$$
|\cR^\dd_j u(t,\theta_t)|\lesssim 2^{-(\alpha+\beta)j}
\|u\|_{\mL^\infty_T(\bC^{\alpha+\beta-\vartheta}_{\dd})}+2^{-(\alpha+\beta)j}\|f\|_{\mL^\infty_T(\bC^\beta_{\dd})},\ j\geq j_0,
$$
and
$$
|\cR^x_j u(t,\theta_t)|\lesssim 2^{-\frac{\gamma+\alpha}{1+\alpha}j}\Big(\|u\|_{\mL^\infty_T(\bC^{(\gamma+\alpha-\vartheta)/(1+\alpha)}_x)}
+\|u\|_{\mL^\infty_T(\bC^{(\gamma-\beta)/(1+\alpha)}_x\bC^{\alpha+\eps}_\v)}
+\|f\|_{\mL^\infty_T(\bC^{\gamma/(1+\alpha)}_x)}\Big),\ j\geq j_0.
$$
Moreover, it is clear that
$$
|\cR^\dd_j u(t,\theta_t)|\leq \|u\|_{\mL^\infty_T},\  j\in\mN_0.
$$
By \eqref{Sur}, Theorem \ref{Th22} and \eqref{In}, for any $\eps'>0$, the above estimates lead to 
\begin{align*}
\|u(t)\|_{\bC^{\alpha+\beta}_{\dd}}
\lesssim\|u(t)\|_{\bB^{\alpha+\beta}_{\dd,\infty}}=
\sup_{j\in\mN_0}2^{(\alpha+\beta)j}\|\cR^\dd_j u(t)\|_\infty
&\lesssim\eps'\|u\|_{\mL^\infty_T(\bC^{\alpha+\beta}_{\dd})}+\|u\|_{\mL^\infty_T}
+\|f\|_{\mL^\infty_T(\bC^\beta_{\dd})},
\end{align*}
and
\begin{align*}
\|u(t)\|_{\bC^{(\gamma+\alpha)/(1+\alpha)}_{x}}&\lesssim \|u(t)\|_{\bB^{(\gamma+\alpha)/(1+\alpha)}_{x,\infty}}
=\sup_{j\in\mN_0}2^{(\gamma+\alpha)j}\|\cR^x_j u(t)\|_\infty\\
&\lesssim\eps'\|u\|_{\mL^\infty_T(\bC^{(\gamma+\alpha)/(1+\alpha)}_x)}+
\|u\|_{\mL^\infty_T(\bC^{(\gamma-\beta)/(1+\alpha)}_x\bC^{\alpha+\eps}_\v)}
+\|f\|_{\mL^\infty_T(\bC^{\gamma/(1+\alpha)}_x)},
\end{align*}
which in turn yield \eqref{Sch1} and \eqref{FW1} by taking $\eps'=1/2$ and \eqref{NM4}.

(ii) For $j\geq 1$, we have
$$
\p_t \cR^x_j u=\sL^{(\alpha)}_{\kappa;\v}\cR^x_j u+b\cdot\nabla \cR^x_j u+\cR^x_j f+[\cR^x_j,\sL^{(\alpha)}_{\kappa;\v}]u+[\cR^x_j,b\cdot\nabla ]u.
$$
For $\theta\in(0,1]$ being small enough, by Schauder's estimate \eqref{Sch1}
and Corollary \ref{Cor34}, we have
\begin{align}\label{JH8}
\begin{split}
\|\cR^x_j u\|_{\mL^\infty_T(\bC^{\alpha+\theta\beta}_{\dd})}
&\lesssim \|\cR^x_j f\|_{\mL^\infty_T(\bC^{\theta\beta}_{\dd})}
+\|[\cR^x_j,\sL^{(\alpha)}_{\kappa;\v}]u\|_{\mL^\infty_T(\bC^{\theta\beta}_{\dd})}
+\|[\cR^x_j,b\cdot\nabla ]u\|_{\mL^\infty_T(\bC^{\theta\beta}_{\dd})}\\
&\lesssim \|\cR^x_j f\|_{\mL^\infty_T(\bC^{\theta\beta}_{\dd})}+2^{-\frac{(1-\theta)\gamma}{1+\alpha}j}
\Big(\eps\|u\|_{\mL^\infty_T(\bC^{(\alpha+(1-\theta)\gamma+\theta\beta)/(1+\alpha)}_x)}+\|u\|_{\mL^\infty_T(\bC^{\alpha+\beta}_\v)}\Big),
\end{split}
\end{align}
where $j\geq j_0$ and $\eps\in(0,1)$.
On the other hand, for $j=0,\cdots,j_0$, 
$$
\|\cR^x_j u\|_{\mL^\infty_T(\bC^{\alpha+\theta\beta}_{\dd})}\lesssim\|\cR^x_j u\|_{\mL^\infty_T(\bC^{(\alpha+\theta\beta)/(1+\alpha)}_x)}
+\|\cR^x_j u\|_{\mL^\infty_T(\bC^{\alpha+\theta\beta}_\v)}\lesssim\|u\|_{\mL^\infty_T(\bC^{\alpha+\beta}_\v)},
$$
and also,
\begin{align}\label{JH9}
\begin{split}
\sup_{j\in\mN_0} 2^{\frac{(1-\theta)\gamma}{1+\alpha}j}\|\cR^x_j f\|_{\mL^\infty_T(\bC^{\theta\beta}_{\dd})}
&\stackrel{\eqref{HQ7}}{\lesssim}\|f\|_{\mL^\infty_T(\bC^{(\theta\beta+(1-\theta)\gamma)/(1+\alpha)}_x)}
+\|f\|_{\mL^\infty_T(\bC^{(1-\theta)\gamma/(1+\alpha)}_x\bC^{\theta\beta}_\v)}\\
&\stackrel{\eqref{GA66}}{\lesssim}\|f\|_{\mL^\infty_T(\bC^{\gamma/(1+\alpha)}_x)}+\|f\|_{\mL^\infty_T(\bC^{\beta}_\v)}.
\end{split}
\end{align}
Hence, by \eqref{HQ7}, \eqref{JH8}, \eqref{JH9} and \eqref{Sch1}, we obtain that for any $\eps\in(0,1)$,
\begin{align*}
&\|u\|_{\mL^\infty_T(\bC^{(\alpha+(1-\theta)\gamma+\theta\beta)/(1+\alpha)}_x)}
+\|u\|_{\mL^\infty_T(\bC^{(1-\theta)\gamma/(1+\alpha)}_x\bC^{\alpha+\theta\beta}_\v)}\lesssim
\sup_{j\in\mN_0} 2^{\frac{(1-\theta)\gamma}{1+\alpha}j}\|\cR^x_j u\|_{\mL^\infty_T(\bC^{\alpha+\theta\beta}_{\dd})}\\
&\qquad\lesssim \|f\|_{\mL^\infty_T(\bC^{\gamma/(1+\alpha)}_x)}+\|f\|_{\mL^\infty_T(\bC^{\beta}_\v)}
+\eps\|u\|_{\mL^\infty_T(\bC^{(\alpha+(1-\theta)\gamma+\theta\beta)/(1+\alpha)}_x)},
\end{align*}
which implies by taking $\eps$ small enough,
$$
\|u\|_{\mL^\infty_T(\bC^{(1-\theta)\gamma/(1+\alpha)}_x\bC^{\alpha+\theta\beta}_\v)}\lesssim
\|f\|_{\mL^\infty_T(\bC^{\gamma/(1+\alpha)}_x)}+\|f\|_{\mL^\infty_T(\bC^{\beta}_\v)}.
$$
Substituting this into \eqref{FW1} with $\theta$ being small enough, we obtain the desired estimate \eqref{Sch}.
\end{proof}

\br\label{Re610}
The restriction of $\alpha\in(1,2)$ is only used in Lemma \ref{Le54}, which is caused by the moment problem due to $1+\vartheta<\alpha$. 
In particular, if $b^{(1)}(t,x,\v)=\v+b^{(1)}(t,x)$, then under the following restrictions:
\begin{align}\label{DS2}
\tfrac{1+\vartheta}{1+\alpha}<\alpha,\ \tfrac{\gamma}{1+\alpha}<\alpha,\ \alpha+\beta>1,
\end{align}
which implies $\alpha>\frac{\sqrt{5}-1}{2}$, Theorem \ref{Th66} still holds for $\alpha\in(\frac{\sqrt{5}-1}{2},1]$. 
Here we conjecture that the moment restriction is superfluous.
At this moment we do not know how to drop it. Such a problem also appears in \cite{Cha-Me-Pr}.
Moreover, if $b(t,x,\v)=(\v, 0)$, which corresponds to the kinetic equation \eqref{Kinetic}, 
then Theorem \ref{Th66} holds for all $\alpha\in(0,2)$.
\er

We have the following existence of classical solutions.
\bt\label{Th603}
Let $\alpha\in(1,2)$ and $\beta\in(0,1), \vartheta\in(0,\alpha-1)$, $\gamma\in(1,1+\alpha)$.
Under {\bf (H$^{\alpha,\vartheta}_{\beta,\gamma}$)}, for any $f\in\mL^\infty_{loc}(\bC^{\gamma/(1+\alpha)}_x\cap\bC^{\beta}_\v)$,
there is a unique classical solution $u$ in the sense of Definition \ref{Def1}
 such that for any $T>0$ and some $C>0$ being independent of $\lambda>0$,
\begin{align}\label{Sch0}
\|u\|_{\mL^\infty_T(\bC^{(\gamma+\alpha)/(1+\alpha)}_x\cap\bC^{\alpha+\beta}_\v)}
\leq C\|f\|_{\mL^\infty_T(\bC^{\gamma/(1+\alpha)}_x\cap \bC^{\beta}_\v)},\
\|u\|_{\mL^\infty_T}\leq \lambda^{-1}\|f\|_{\mL^\infty_T}.
\end{align}
\et
\begin{proof}
Let $(\rho_n)_{n\in\mN}$ and $(\rho'_n)_{n\in\mN}$ be the usual mollifiers in $\mR^{3d}$ and $\mR^{2d}$ respectively. 
Define
$$
\kappa_n(t,x,\v,w):=\kappa(t,\cdot)*\rho_n(x,\v,w),\ n\in\mN,
$$ 
and
$$
b_n(t,x,\v):=b(t,\cdot)*\rho'_n(x,\v),\ \ f_n(t,x,\v):=f(t,\cdot)*\rho'_n(x,\v),\ n\in\mN.
$$
Fix $\lambda, T>0$. By Theorem \ref{Th82} in appendix, there is a unique smooth
$u_n\in C([0,T]; \sC^2(\mR^{2d}))$
solving the following PDE:
\begin{align}\label{APP}
\p_t u_n=\sL^{(\alpha)}_{\kappa_n;\v} u_n+b_n\cdot\nabla u_n-\lambda u_n+f_n,\  \ u_n(0)=0.
\end{align}
Under {\bf (H$^{\alpha,\vartheta}_{\beta,\gamma}$)}, by Theorem \ref{Th66},
there is a constant $C>0$ such that for  all $n\in\mN$,
\begin{align}\label{KA1}
\|u_n\|_{\mL^\infty_T(\bC^{(\gamma+\alpha)/(1+\alpha)}_x\cap \bC^{\alpha+\beta}_\v)}
\leq C\|f_n\|_{\mL^\infty_T(\bC^{\gamma/(1+\alpha)}_x\cap\bC^{\beta}_\v)}
\leq C\|f\|_{\mL^\infty_T(\bC^{\gamma/(1+\alpha)}_x\cap\bC^{\beta}_\v)}.
\end{align}
Moreover, since $\alpha\in(1,2)$ and $\gamma\in(1,1+\vartheta)$, we also have for some $\eps>0$,
$$
\|\nabla u_n\|_{\mL^\infty_T(\bC^{\eps})}\leq C
\|u_n\|_{\mL^\infty_T(\bC^{(\gamma+\alpha)/(1+\alpha)}_x\cap \bC^{\alpha+\beta}_\v)}\leq C\|f\|_{\mL^\infty_T(\bC^{\gamma/(1+\alpha)}_x\cap\bC^{\beta}_\v)}.
$$
Hence, from approximation equation \eqref{APP} and the above uniform estimates, one sees that
$$
\sup_n\|\p_t u_n\cdot \1_{|x|+|\v|\leq m}\|_{\mL^\infty_T}\leq C_m,\ \ m\in\mN.
$$
Thus by Ascolli-Arzela's theorem and a standard diagonalization argument, there are subsequence $n_k$ 
and continuous function $u:[0,T]\times\mR^{2d}\to\mR$ such that for each $m\in\mN$,
$$
\lim_{k\to\infty}\sup_{t\in[0,T],|x|+|\v|\leq m}|u_{n_k}(t,x,\v)-u(t,x,\v)|=0.
$$
Moreover, we also have
$$
u\in \mL^\infty_T(\bC^{(\gamma+\alpha)/(1+\alpha)}_x\cap \bC^{\alpha+\beta}_\v).
$$ 
In fact, by \eqref{LJ1} and Fatou's lemma, we have 
\begin{align*}
\|u(t)\|_{\bC^{\alpha+\beta}_\v}&\lesssim\sup_{j\geq 0} 2^{(\alpha+\beta)j}\|\cR^\v_j u(t)\|_\infty
\lesssim\sup_{j\geq 0} 2^{(\alpha+\beta)j}\varliminf_{n\to\infty}\|\cR^\v_j u_n(t)\|_\infty\\
&\lesssim \varliminf_{n\to\infty}\|u_n(t)\|_{\bC^{\alpha+\beta}_\v}\stackrel{\eqref{KA1}}{\lesssim}
\|f\|_{\mL^\infty_T(\bC^{\gamma/(1+\alpha)}_x\cap\bC^{\beta}_\v)}.
\end{align*}
Noticing the following interpolation inequality (see \cite[Theorem 3.2.1]{Kr}),
$$
\|\nabla f\|_\infty\leq C\|f\|^{1/(1+\eps)}_{\bC^{1+\eps}}\|f\|^{\eps/(1+\eps)}_\infty,
$$
we further have
$$
\lim_{k\to\infty}\sup_{t\in[0,T],|x|+|\v|\leq m}|\nabla u_{n_k}(t,x,\v)-\nabla u(t,x,\v)|=0.
$$
By taking limits for equation \eqref{APP}, we obtain that $u$ satisfies \eqref{PDE0} in the sense of Definition \ref{Def1}.
By \eqref{NM4}, we complete the proof.
\end{proof}

\br
If we do not assume $\gamma>1$ in Theorem \ref{Th603}, 
then under {\bf (H$^{\alpha,\vartheta}_{\beta,\gamma}$)}, for any $f\in \mL^\infty_T(\bC^{\beta}_a)$, 
we can show the existence of $u\in\mL^\infty_T(\bC^{\alpha+\beta}_a)$ solving PDE \eqref{PDE0} in the distributional sense
since $b\cdot\nabla u$ is a distribution under the above regularity. 
\er

\section{Degenerate SDEs with H\"older drifts}

\subsection{Pathwise uniqueness of SDEs with multiplicative L\'evy noises}
Let $L^{(\alpha)}_t$ be a symmetric and rationally invariant $\alpha$-stable process with $\alpha\in(1,2)$
 on some probability space $(\Omega,\sF,\mP)$, so that whose generator is given by the fractional Laplacian $\Delta^{\alpha/2}$.
In this section we consider the following degenerate SDE with jumps in $\mR^{2d}$:
\begin{align}\label{SDE3}
\dif Z_{s,t}=b(t,Z_{s,t})\dif t+(0,\sigma(t,Z_{s,t})\dif L^{(\alpha)}_t),\ \ Z_{s,s}=z\in\mR^{2d},\ t\geq s\geq 0,
\end{align}
where $\sigma:\mR_+\times\mR^{2d}\to \mR^d\otimes\mR^d$ 
and $b:\mR_+\times\mR^{2d}\to\mR^{2d}$ are measurable functions satisfying
\begin{enumerate}[{\bf (${\widetilde {\bf H}}^{\alpha,\vartheta}_{\beta,\gamma}$)}]
\item $\sigma$ is Lipschitz continuous in $x$ uniformly in $t$, and for some $c_0\geq 1$ and all $t\geq 0$,
$$
c_0^{-1}|\xi|\leq |\sigma(t,z)\xi|\leq c_0|\xi|,\ \ \xi\in\mR^d,\ z\in\mR^{2d},
$$
and for some $\vartheta\in(0,\alpha-1)$, $\gamma\in(1,1+\alpha)$ and $\beta\in(0,1)$, 
$$
|b(t,0)|+[b(t,\cdot)]_{\bC^{\gamma/(1+\alpha)}_x}+\|\nabla_\v b^{(1)}(t,\cdot)\|_{\bC^{\vartheta}_\v}+[b^{(2)}(t,\cdot)]_{\bC^{\beta}_\v}\leq c_0.
$$
Moreover, \eqref{GA4} holds.
\end{enumerate}

Let $N(\dif t,\dif w)$ be the Poisson random measure associated with ${L^{(\alpha)}}$, i.e.,
$$
N((0,t]\times\Gamma):=\sum_{0<s\leq t}1_\Gamma(L^{(\alpha)}_s-L^{(\alpha)}_{s-}), \quad t>0, \
\Gamma\in \sB(\mR^d\setminus\{0\}).
$$
Let $\tilde N(\dif t,\dif w):=N(\dif t,\dif w)-\dif t\dif w/|w|^{d+\alpha}$ be the compensated Poisson random measure.
By the L\'evy-It\^o decomposition, we can write for each $t>0$,
$$
L^{(\alpha)}_t=\int^t_0\!\!\!\int_{\mR^d}w\tilde N(\dif s,\dif w).
$$
Thus, if we let $Z_{s,t}=(X_{s,t},V_{s,t})$, then SDE \eqref{SDE3} can be written as
$$
\left\{
\begin{aligned}
&\dif X_{s,t}=b^{(1)}(t,Z_{s,t})\dif t,\ \ (X_{s,s},V_{s,s})=(x,\v),\\
&\dif V_{s,t}=b^{(2)}(t,Z_{s,t})\dif t+\int_{\mR^d}\sigma(t,Z_{s,t})w\tilde N(\dif t,\dif w).
\end{aligned}
\right.
$$
In particular, the generator of this SDE is given by $\sL^{(\alpha)}_{\sigma;\v}+b\cdot\nabla$ with
\begin{align*}
\sL^{(\alpha)}_{\sigma;\v} f(x,\v)&={\rm p.v.}\int_{\mR^d}\Big(f(x,\v+\sigma(t,z)w)-f(x,\v)\Big)\frac{\dif w}{|w|^{d+\alpha}}\\
&={\rm p.v.}\int_{\mR^d}\Big(f(x,\v+w)-f(x,\v)\Big)\kappa(t,z,w)\frac{\dif w}{|w|^{d+\alpha}},
\end{align*}
where $z=(x,\v)$ and
$$
\kappa(t,z,w):=\det(\sigma^{-1}(t,z))|w|^{d+\alpha}/|\sigma^{-1}(t,z)w|^{d+\alpha}.
$$
Under {\bf (${\widetilde {\bf H}}^{\alpha,\vartheta}_{\beta,\gamma}$)}, it is easy to see that for some $c_0>1$,
$$
c_0^{-1}\leq \kappa(t,z,w)\leq c_0,\ \ |\kappa(t,z,w)-\kappa(t,z',w)|\leq c_0|z-z'|.
$$

We have the following result.
\bt\label{Th71}
Let $\alpha\in(1,2)$, $\vartheta\in(0,\alpha-1)$, $\gamma\in(1+\frac{\alpha}{2},1+\alpha)$ and $\beta\in(1-\frac{\alpha}{2},1)$. 
Under {\bf (${\widetilde {\bf H}}^{\alpha,\vartheta}_{\beta,\gamma}$)}, for each $s\geq 0$ and $z\in\mR^{2d}$,
there exists a unique strong solution $(Z_{s,t})_{t\geq s}$ to SDE \eqref{SDE3}.
%Moreover, $\{z\mapsto Z_t(z), t\geq 0\}$ forms a $C^1$-stochastic diffeomorephism flow.
\et
\begin{proof}
Since the coefficients are continuous and linear growth, the existence of a 
solution is well-known (cf. \cite{Ja}). By Yamada-Watanabe's theorem (cf.\cite{Ik-Wa}), 
it suffices to show the pathwise uniqueness.
Without loss of generality, we assume $s=0$ and simply write
$$
Z_t:=Z_{0,t}.
$$
Since $b$ is unbounded, to construct Zvonkin's transformation (cf. \cite{Zv}), we need to cutoff $b$.
For $m\in\mN$, let $\chi_m$ be a smooth cutoff function in $\mR^{2d}$ with
$$
\chi_m(z)=1,\ |z|\leq m,\ \chi_m(z)=0,\ |z|>m+1.
$$
Fix $T>0$ and $m\in\mN$. Consider the following backward equation:
\begin{align}\label{EG4}
\p_s\u^m_\lambda+\sL^{(\alpha)}_{\sigma;\v}\u^m_\lambda-\lambda \u^m_\lambda+b\cdot\nabla \u^m_\lambda+b\chi_m=0,\ \u^m_\lambda(T,\cdot)=0.
\end{align}
Under {\bf (${\widetilde {\bf H}}^{\alpha,\vartheta}_{\beta,\gamma}$)}, by Theorem \ref{Th603}, there is a 
unique classical solution $\u^m_\lambda$ to the above equation with regularity:
for some $C_m\geq 1$ and all $\lambda\geq 1$,
$$
\|\u^m_\lambda\|_{\mL^\infty_T(\bC^{(\gamma+\alpha)/(1+\alpha)}_x\cap\bC^{\alpha+\beta}_\v)}
\leq C_m,\quad \|\u^m_\lambda\|_{\mL^\infty_T}\leq C_m\lambda^{-1}.
$$
Since $\gamma\in(1+\frac{\alpha}{2},1+\alpha)$ and $\beta\in(1-\frac{\alpha}{2},1)$,  
by interpolation inequalities \eqref{In} and \eqref{GA66}, there is an $\eps_0>0$ small enough such that 
$$
\|\u^m_\lambda\|_{\mL^\infty_T(\bC^{1+\eps_0}_x\bC^{\alpha/2+\eps_0}_\v)}
+\|\u^m_\lambda\|_{\mL^\infty_T(\bC^{1+\alpha/2+\eps_0}_\v)}\leq c_\lambda,
$$
where $c_\lambda\to 0$ as $\lambda\to\infty$.
In particular, for $\lambda\geq 1$ large enough,
\begin{align}\label{UH7}
\|\nabla \u^m_\lambda\|_{\mL^\infty_T(\bC^{\alpha/2+\eps_0}_\v)}\leq \tfrac{1}{2}.
\end{align}
Define
$$
\Phi^m_t(z):=z+\u^m_\lambda(t,z).
$$
By \eqref{UH7}, one sees that
$$
z\mapsto\Phi^m_t(z) \mbox{ is a diffeomorphism on $\mR^{2d}$}
$$
and
$$
\|\nabla\Phi^m_t\|_\infty+\|\nabla(\Phi^m_t)^{-1}\|_\infty\leq 2.
$$
Moreover, by \eqref{EG4} we have
\begin{align}\label{KH1}
\p_s\Phi^m+\sL^{(\alpha)}_{\sigma;\v}\Phi^m+b\cdot\nabla\Phi^m=b(1-\chi_m)+\lambda\u^m_\lambda.
\end{align}
Let $Z_t$ and $Z'_t$ be two solutions of SDE \eqref{SDE3} defined on the same probability space with the same starting point $z$.  Define a stopping time
$$
\tau_m:=\inf\Big\{t>0: |Z_t|\wedge|Z'_t|\geq m+1\mbox{ or } |\Delta L_t|>m\Big\},
$$
and let
$$
g^m_s(z,w):=\Phi^m_s(z+(0,\sigma(s,z)w))-\Phi^m_s(z).
$$
By  It\^o's formula and \eqref{KH1}, one sees that 
\begin{align*}%\label{SDE4}
\Phi^m_{t\wedge\tau_m}(Z_{t\wedge\tau_m})-\Phi^m_0(z)
&=\int^{t\wedge\tau_m}_0\!\!\int_{\mR^d}
\Big(\Phi^m_s(Z_{s-}+(0,\sigma(s,Z_{s-})w))-\Phi^m_s(Z_{s-})\Big)\tilde N(\dif s,\dif w)\\
&\quad+\int^{t\wedge\tau_m}_0\Big(\p_s\Phi^m+b\cdot\nabla\Phi^m+\sL^{(\alpha)}_{\sigma;\v}\Phi^m\Big)(s,Z_s)\dif s\\
&=\int^{t\wedge\tau_m}_0\!\!\int_{B_m}g^m_s(Z_{s-},w)\tilde N(\dif s,\dif w)+\lambda\int^{t\wedge\tau_m}_0\u^m_\lambda(s,Z_s)\dif s,
\end{align*}
where we have used that $(b(1-\chi_m))(s,Z_s)=0$ for $s<\tau_m$ and
$$
\int^{t\wedge\tau_m}_0\!\!\int_{B_{m}^c}g^m_s(Z_{s-},w)\tilde N(\dif s,\dif w)
=\sum_{0<s\leq t\wedge\tau_m}g^m_s(Z_{s-},\Delta L_s)\cdot \1_{|\Delta L_s|>m}=0.
$$
Noticing that
$$
\nabla_z g^m_s(z,w)=(\nabla\Phi^m_s)(z+(0,\sigma(s,z)w))(\mI+(0,\nabla\sigma(s,z)w))-\nabla\Phi^m_s(z),
$$ 
by \eqref{UH7} we have
$$
\|\nabla_z g^m_s(\cdot,w)\|_\infty\leq 2(\|\sigma\|_\infty|w|)^{\alpha/2+\eps_0}+2\|\nabla_z\sigma\|_\infty|w|.
$$
Hence, by the isometry formula of stochastic integrals, 
\begin{align*}
\mE|Z_{t\wedge\tau_m}&-Z'_{t\wedge\tau_m}|^2\leq 4\mE|\Phi^m_{t\wedge\tau_m}(Z_{t\wedge\tau_m})-\Phi^m_{t\wedge\tau_m}(Z'_{t\wedge\tau_m})|^2\\
&\lesssim\mE\int^{t\wedge\tau_m}_0\!\!\int_{B_m}|g^m_s(Z_{s-},w)-g^m_s(Z'_{s-},w)|^2\frac{\dif w}{|w|^{d+\alpha}}\dif s\\
&+\lambda\mE\int^{t\wedge\tau_m}_0|\u^m_\lambda(s,Z_s)-\u^m_\lambda(s,Z'_s)|^2\dif s\\
&\lesssim \mE\int^{t\wedge\tau_m}_0|Z_s-Z'_s|^2\dif s\left(\int_{B_m}\frac{(|w|^{\alpha+2\eps_0}+|w|^2)\dif w}{|w|^{d+\alpha}}+\lambda\right),
\end{align*}
which yields by Gronwall's inequality
$$
Z_{t\wedge\tau_m}=Z'_{t\wedge\tau_m},\ \ t\geq 0.
$$
Finally, letting $m\to\infty$, we obtain the pathwise uniqueness.
\end{proof}

\br
By suitable localization technique, we can directly construct the solution by Zvonkin's transformation without using Yamada-Watanabe's theorem.
\er

\subsection{$C^1$-stochastic diffeomorphism flows for SDEs with additive L\'evy noise}
In this subsection we consider the $C^1$-stochastic diffeomorphism flows property 
for SDE \eqref{SDE3} with additive L\'evy noises. We introduce the following spaces:
For a Fr\'echet space $\mF$ and time interval $I$, define
$$
C(I;\mF):=\{f: I\to\mF \mbox{ is continuous}\},\ \ 
D(I;\mF):=\{f: I\to\mF \mbox{ is c\`adl\`ag}\}.
$$
For $k\in\mN_0$, let $\mC^k$ be the Fr\'echet space of all 
$k$-order continuous differentiable functions with Fr\'echet metric:
$$
d(f,g):=\sum_{j=0}^k\sum_{n\in\mN} 2^{-n}\left(1\wedge\sup_{|x|<n}|\nabla^j f(x)-\nabla^j g(x)|\right).
$$
We have
\bt\label{Th72}
Let $\alpha\in(1,2)$, $\gamma\in(1+\frac{\alpha}{2},1+\alpha)$ and $\beta\in(1-\frac{\alpha}{2},1)$. Assume 
$\sigma\equiv 1$ and 
\begin{align}\label{DS1}
b(t,x,\v)=(\v+b^{(1)}(t,x), b^{(2)}(t,x,\v)),
\end{align}
with
$$
b^{(1)}\in\mL^\infty_{loc}(\bC^{\gamma/(1+\alpha)}_x),\ \ b^{(2)}\in\mL^\infty_{loc}(\bC^{\gamma/(1+\alpha)}_x\cap\bC^\beta_\v).
$$
Then the unique strong solution $\{Z_{s,t}(z), t>s\geq 0, z\in\mR^{2d}\}$ of SDE \eqref{SDE3} forms a $C^1$-stochastic diffeomorphism flow. More precisely,
there is a null set $\cN$ such that for all $\omega\notin\cN$,
\begin{enumerate}[(i)]
\item For all $0\leq s<r<t$, it holds that
$$
Z_{s,t}(z,\omega)=Z_{r,t}(Z_{s,r}(z,\omega),\omega),\ \ \forall z\in\mR^{2d},
$$
and
$$
z\mapsto Z_{s,t}(z,\omega)\mbox{ is a $C^1$-diffeomorphism on $\mR^{2d}$.}
$$
\item %$s\mapsto Z_{s,t}(\cdot,\omega)$ and 
$t\mapsto Z_{s,t}(\cdot,\omega)\in D([s,\infty); \mC^1)$ and $s\mapsto Z_{s,t}(\cdot,\omega)\in D([0,t]; \mC^1)$.
\end{enumerate}
\et
\begin{proof}
Fix $T,\lambda>0$. Consider the following backward equation:
$$
\p_s\u_\lambda+\Delta^{\alpha/2}_{\v}\u_\lambda-\lambda \u_\lambda+b\cdot\nabla \u_\lambda+b=0,\ \u_\lambda(T,\cdot)=0,
$$
where $\Delta^{\alpha/2}_{\v}$ is the fractional Laplacian acting on the variable $\v$.
Under the assumptions of the theorem, by Theorem \ref{Th603}, there is a 
unique classical solution $\u^m_\lambda$ to the above equation with regularity:
for some $C\geq 1$ and all $\lambda\geq 1$,
$$
\|\u_\lambda\|_{\mL^\infty_T(\bC^{(\gamma+\alpha)/(1+\alpha)}_x\cap\bC^{\alpha+\beta}_\v)}
\leq C,\quad \|\u_\lambda\|_{\mL^\infty_T}\leq C\lambda^{-1}.
$$
As in Theorem \ref{Th71}, by \eqref{In} and \eqref{GA66}, 
there are $\lambda\geq 1$ large enough and $\eps_0>0$ such that
\begin{align}\label{UH77}
\|\nabla \u_\lambda(t,\cdot,\cdot)\|_{\bC^{\alpha/2+\eps_0}_\v}\leq \tfrac{1}{2},\ \|\u_\lambda\|_{\mL^\infty_T(\bC^{1+\eps_0}_x\bC^{\alpha/2+\eps_0}_\v)}
+\|\u_\lambda\|_{\mL^\infty_T(\bC^{1+\alpha/2+\eps_0}_\v)}<\infty.
\end{align}
Define
$$
\Phi_t(z):=z+\u_\lambda(t,z).
$$
By \eqref{UH77}, one sees that
$$
z\mapsto\Phi_t(z) \mbox{ is a diffeomorphism on $\mR^{2d}$}
$$
and
$$
\|\nabla\Phi_t\|_\infty+\|\nabla\Phi_t^{-1}\|_\infty\leq 2.
$$
Let
$$
\widetilde Z_t:=\Phi_t(Z_t).
$$
By It\^o's formula, one sees that 
\begin{align}\label{SDE9}
\widetilde Z_t&=\Phi_t(Z_t)=\Phi_0(z)+\int^t_0\Big(\p_s\Phi+b\cdot\nabla\Phi+\Delta^{\alpha/2}_{\v}
\Phi\Big)(s,Z_s)\dif s\no\\
&+\int^t_0\!\!\int_{\mR^d}\Big(\Phi_s(Z_{s-}+(0,w))-\Phi_s(Z_{s-})\Big)\tilde N(\dif s,\dif w)\no\\
&=\tilde z+\int^t_0\widetilde b(s,\widetilde Z_s)\dif s+\int^t_0\!\!\!\int_{\mR^d}\widetilde g_s(\widetilde Z_{s-},w)\tilde N(\dif s,\dif w),
\end{align}
where $\tilde z:=\Phi_0(z)$ and
$$
\widetilde b(t,z):=\lambda \u_\lambda(t,\Phi^{-1}_t(z)),\quad \widetilde g_t(z,w):=\Phi_t(\Phi^{-1}_t(z)+(0,w))-z.
$$
In particular, $\{Z_{s,t}, t\geq 0\}$ solves SDE \eqref{SDE3} if and only if $\{\widetilde Z_{s,t}, t\geq 0\}$
solves SDE \eqref{SDE9} (see \cite[Lemma 3.4]{Ch-So-Zh}).
\\
\\
{\it Claim:} There are $\eps_1,\eps_2\in(0,\eps_0)$ and constant $C>0$ such that
for all $t\in[0,T]$ and $w\in\mR^d$,
$$
\|\nabla_z \widetilde g_t(\cdot,w)\|_{\bC^{\eps_1}}\leq C(|w|^{\alpha/2+\eps_2}\wedge 1),\ \ 
\|\nabla_z \widetilde b\|_{\mL^\infty_T(\bC^{\eps_1})}<\infty.
$$
Indeed, noticing that
\begin{align*}
&\nabla_z\widetilde g_t(z,w)=(\nabla\Phi_t)(\Phi^{-1}_t(z)+(0,w))\cdot\nabla\Phi^{-1}_t(z)-\mI\\
&\quad=\Big((\nabla\Phi_t)(\Phi^{-1}_t(z)+(0,w))-\nabla\Phi_t\circ\Phi^{-1}_t(z)\Big)\cdot\nabla\Phi^{-1}_t(z),
\end{align*}
by \eqref{UH77} we have
$$
\|\nabla_z\widetilde  g_t(\cdot,w)\|_\infty\leq \|\nabla\Phi_t\|_{\bC^{\alpha/2+\eps_0}_\v}(|w|^{\alpha/2+\eps_0}\wedge 1)
\|\nabla\Phi^{-1}_t\|_\infty\lesssim 1\wedge|w|^{\alpha/2+\eps_0},
$$
and by definition,
\begin{align*}
[\nabla_z\widetilde  g_t(\cdot,w)]_{\bC^{\eps_0}}
&\leq [\nabla\Phi_t(\Phi^{-1}_t(\cdot)+(0,w))]_{\bC^{\eps_0}}
\|\nabla\Phi^{-1}_t\|_\infty+\|\nabla\Phi_t\|_\infty[\nabla\Phi^{-1}_t]_{\bC^{\eps_0}}
\lesssim 1.
\end{align*}
Thus we obtain the first claim by standard interpolation technique. The second one is easy by definition and \eqref{UH77}.
\\
\\
By the above claim and \cite[Theorem 4.1]{Fu-Ku}, the unique solutions of SDE \eqref{SDE9} and so 
SDE \eqref{SDE3} define a $C^1$-stochastic diffeomorphism flow
and (i), (ii) hold. See also \cites{Pr12, Ch-So-Zh} for more details.
\end{proof}

\br
Although we here only consider the symmetric and rotationally invariant $\alpha$-stable additive noise, 
by the same argument as used in \cite{Ch-So-Zh}, it is possible to consider more general additive L\'evy noises, even cylindrical L\'evy noises.
On the other hand, by Remark \ref{Re610},  the $\alpha$ in Theorem \ref{Th72} in fact can be less than $1$, but with a lower bound $(\sqrt{17}-1)/4\approx 0.78078$. Indeed, by
restriction \eqref{DS2} and $\gamma\in(1+\frac{\alpha}{2},1+\alpha)$, we need to require 
$$
\alpha^2+\alpha>1+\tfrac{\alpha}{2}\Rightarrow\alpha>(\sqrt{17}-1)/4.
$$
\er

\subsection{Application to random transport equations with H\"older coefficients}
In this subsection we apply Theorem \ref{Th72} to a random transport equation with H\"older coefficient.
First of all, the following corollary is an easy consequence of Theorem \ref{Th72}.
\bc\label{Cor73}
Let $\alpha\in(1,2)$ and $\gamma\in(\frac{2+\alpha}{2(1+\alpha)},1)$. Assume that
$$
b^{(1)}:\mR_+\times\mR^d\to\mR^d\in \mL^\infty_{loc}(\bC^\gamma).
$$
Consider the following random ODE
\begin{align}\label{ODE0}
\dif Y_{s,t}(x,\omega)/\dif t=b^{(1)}(t,Y_{s,t}(x,\omega))+L^{(\alpha)}_t(\omega),\ t>s, \ Y_{s,s}=x. 
\end{align}
For $\mP$-almost all $\omega$, there exists a family of solutions $\{Y_{s,t}(x,\omega), x\in\mR^d, 0\leq s<t<\infty\}$ to the above random ODE
so that 
\begin{enumerate}[(i)]
\item For each $s<r<t$, it holds that
\begin{align}\label{DP1}
Y_{s,t}(x,\omega)=Y_{r,t}(Y_{s,r}(x,\omega),\omega),\ \ x\in\mR^d,
\end{align}
and
\begin{align}\label{DP2}
x\mapsto Y_{s,t}(x,\omega)\mbox{ is a $C^1$-diffeomorphism on $\mR^d$.}
\end{align}
\item For each $s\geq 0$, $t\mapsto Y_{s,t}(\cdot,\omega)\in C([s,\infty),\mC^0)\cap D([s,\infty),\mC^1)$.
\end{enumerate}
\ec
\begin{proof}
For $z=(x,\v)$, by Theorem \ref{Th72} let 
$$
Z_{s,t}(z):=(X_{s,t}(x,\v), L_t-L_s+\v)
$$ 
solve the degenerate SDE \eqref{SDE3} with 
$$
\sigma\equiv 1,\ \ b(t,x,\v):=(\v+b^{(1)}(t,x),0).
$$
Now we define
$$
Y_{s,t}(x):=X_{s,t}(x, L_s)
$$
Clearly, it solves ODE \eqref{ODE0}. Now we check (i) and (ii) hold.
By (i) of Theorem \ref{Th72}, we have
\begin{align*}
(Y_{s,t}(x), L_t)&=(X_{s,t}(x,L_s), L_t)=Z_{s,t}(x, L_s)=Z_{r,t}\circ Z_{s,r}(x, L_s)\\
&=Z_{r,t}(X_{s,r}(x,L_s), L_r)=(X_{r,t}(X_{s,r}(x,L_s), L_r), L_t)
=(Y_{r,t}\circ Y_{s,r}(x), L_t),
\end{align*}
which implies \eqref{DP1} and \eqref{DP2} by (i) of Theorem \ref{Th72}.
Finally, (ii) follows by (ii) of Theorem \ref{Th72} and equation \eqref{ODE0}.
\end{proof}
\br
Here an open question is to show Davie's uniqueness \cite{Da} for the above random ODE, that is, for almost all $\omega$, ODE \eqref{ODE0}
has a unique solution. See \cite{Pr18} for the study of random ODE $\dif Y_{t}/\dif t=b(t,Y_{t}+L^{(\alpha)}_t)$. We will study this problem in a future work.
\er

We need the following real analysis result, which can be proven by the completely same method as in \cite[p149, Theorem 7.21]{Ru}. We omit the details.
\bl\label{Le77}
Let $f:[a,b]\to \mR$ be a continuous function. Assume that for each point $t\in[a,b)$, the right derivative $f'_+:=\lim_{\eps\downarrow 0}(f(t+\eps)-f(t))/\eps$ exists and 
$f'_+\in L^1([a,b])$. Then $f$ is absolutely continuous on $[a,b]$.
\el

Now we can state the following result.
\bt\label{Th75}
Let $\alpha\in(1,2)$. Assume 
$b\in C(\mR_+\times\mR^d)\cap\mL^\infty_{loc}(\bC^\gamma)$ with $\gamma\in(\frac{2+\alpha}{2(1+\alpha)},1)$.
For any $\varphi\in \sC^1$ and almost all $\omega$, there is a unique function $(t,x)\mapsto u(t,x,\omega)\in 
C(\mR_+; \mC^0)\cap D(\mR_+; \mC^1)$ 
so that for each $x\in\mR^d$, $t\mapsto u(t,x,\omega)$ is absolutely continuous and
\begin{align}\label{Tran}
\p_tu(t,x,\omega)+(b(t,x)+L_t(\omega))\cdot\nabla_x u(t,x,\omega)=0,\ \ u(0,x)=\varphi(x).
\end{align}
\et
\begin{proof}
We only show the existence since the uniqueness is standard by the maximum principle (see  \cite[Theorem 6.1]{Ch-HXZ}).
By Corollary \ref{Cor73}, for each $x\in\mR^d$, let $X_{s,t}(x,\omega)$ solve the following random ODE:
\begin{align}\label{EQ1}
X_{s,t}(x,\omega)=x+\int^t_s\Big(b(r,X_{s,r}(x,\omega))+L^{(\alpha)}_r(\omega)\Big)\dif r,\ \ 0\leq s\leq t.
\end{align}
%Let $\sD_{x}$ be the continuous points of $s\mapsto X_{s,T}(x),\nabla_x X_{s,T}(x)$.
%It is easy to see that
%$$
%s\mapsto X_{s,t}(x)\mbox{ is continuous.}
%$$
%In fact, for $s<s'<t$, by equation \eqref{EQ1} we have
%\begin{align*}
%|X_{s,t}-X_{s',t}|&\leq \int^t_{s'}|b(r,X_{s,r})-b(r,X_{s',r})|\dif r+(\|b\|_\infty+\|L^{(\alpha)}\|_\infty)(s'-s)
%\end{align*}
Define
$$
u(t,x,\omega):=\varphi(X^{-1}_{0,t}(x,\omega)).
$$
Clearly, by (ii) of Corollary \ref{Cor73}, we have for almost all $\omega$,
\begin{align}\label{GG4}
u(\cdot,\cdot,\omega)\in C(\mR_+; \mC^0)\cap D(\mR_+; \mC^1),
\end{align}
and by (i) of Corollary \ref{Cor73}, for $\eps>0$,
$$
u(t+\eps,x)=\varphi(X^{-1}_{0,t}\circ X^{-1}_{t,t+\eps}(x))=u(t,X^{-1}_{t,t+\eps}(x)).
$$
Here and below we drop the $\omega$.
Hence,
\begin{align}\label{GK1}
\frac{u(t+\eps,x)-u(t,x)}{\eps}=\frac{X^{-1}_{t,t+\eps}(x)-x}{\eps}\int^1_0(\nabla_x u)\Big(t,\theta X^{-1}_{t,t+\eps}(x)+(1-\theta)x\Big)\dif\theta.
\end{align}
Since $\eps\mapsto X_{t,t+\eps}(\cdot)\in D(\mR_+;\mC^1)$
and $\nabla X^{-1}_{t,t+\eps}(x)=(\nabla X_{t,t+\eps})^{-1}\circ X^{-1}_{t,t+\eps}(x)$, we have
\begin{align}\label{GF0}
\lim_{\eps\downarrow 0}\sup_{|x|\leq R}|\nabla X^{-1}_{t,t+\eps}(x)-\mI|=0,\ \ \forall R>0.
\end{align}
Noticing that
$$
X^{-1}_{t,t+\eps}(x)-x=(x-X_{t,t+\eps}(x))
\cdot\int^1_0\nabla X^{-1}_{t,t+\eps}(\theta x+(1-\theta)X_{t,t+\eps}(x))\dif\theta,
$$
since $(t,x)\mapsto b(t,x)$ is continuous and $t\mapsto L^{(\alpha)}_t$ is right continuous,
by \eqref{EQ1} and \eqref{GF0}, for each $(t,x)\in\mR_+\times\mR^d$, we have
$$
\lim_{\eps\downarrow 0}(X^{-1}_{t,t+\eps}(x)-x)/\eps=-b(t,x)-L^{(\alpha)}_t.
$$
Therefore, by \eqref{GK1} and the continuity of $x\mapsto \nabla_x u(t,x)$,
$$
\p^+_t u(t,x):=
\lim_{\eps\downarrow 0}\frac{u(t+\eps,x)-u(t,x)}{\eps}=-(b(t,x)+L^{(\alpha)}_t)\cdot\nabla_x u(t,x), \ \forall (t,x)\in\mR_+\times\mR^d,
$$
where $\p^+_t u$ stands for the right derivative.
Since $t\mapsto (b(t,x)+L^{(\alpha)}_t)\cdot\nabla_x u(t,x)$ is bounded, by Lemma \ref{Le77},
$t\mapsto u(t,x)$ is absolutely continuous. The proof is complete.
\end{proof}
\br
If $\varphi\in L^\infty(\mR^d)$ and $\div b\in L^1_{loc}(\mR_+\times\mR^d)$, then as in \cite[Theorem 20]{Fl-Gu-Pr}, we can show that $u(t,x):=\varphi(X^{-1}_{0,t}(x))$
is the unique bounded weak solution of transport equation \eqref{Tran}.
\er
\section{Appendix}

In this section we use a probabilistic method to show the existence of a smooth solution to PDE \eqref{PDE0} when the coefficients are smooth.
We first recall the following result proved in \cite{Ha-Pe-Zh}.

\bl\label{Cor72}
Given $d_0\in\mN$ and $c_0>1$, let $\kappa(x,z):\mR^{d_0}\times B_1\to[c_0^{-1},c_0]$ be a smooth function
with bounded derivatives of all orders. 
For any $\alpha\in(0,2)$, there is a measurable map $\Phi(x,z): \mR^{d_0}\times B_1\to B_1$ such that for any nonnegative measurable function $f$,
$$
\int_{B_1}f\circ\Phi(x,z)\frac{\dif z}{|z|^{d+\alpha}}=\int_{B_1}f(z)\kappa(x,z)\frac{\dif z}{|z|^{d+\alpha}}.
$$
Moreover, $\Phi$ enjoys the following properties:
\begin{enumerate}[(i)]
\item $\Phi(x,0)=0$ and if $\kappa(x,-z)=\kappa(x,z)$, then $\Phi(x,-z)=-\Phi(x,z)$.
\item For all $i,j\in\mN_0$, there is a $C_{ij}>0$ such that for all $x\in\mR^{d_0}$ and $z\in B_1$,
$$
|\nabla_x^i\nabla_z^j\Phi(x,z)|\leq C_{ij}|z|^{1-j},
$$
where $C_{ij}$ is a polynomial of $\|\nabla^m_x\nabla^n_z\kappa\|_\infty$, $m=1,\cdots$, $i, n=0,\cdots,j$.
\end{enumerate}
\el

\bt\label{Th82}
Suppose that $\kappa$ and $b$ satisfy that for any $m\in\mN$ and $t>0$,
$$
c_0^{-1}\leq \kappa(t,x,\v,w)\leq c_0,\ \ \|\nabla^m\kappa(t,\cdot)\|_\infty+\|\nabla^mb(t,\cdot)\|_\infty\leq c_m.
$$
Then there is a classical solution $u\in\cap_{m\in\mN}C(\mR_+; \sC^m)$ to PDE \eqref{PDE0}.
\et
\begin{proof}
We decompose the operator $\sL^{(\alpha)}_{\kappa;\v}$ as two parts: small jumps $\widetilde{\sL^{(\alpha)}_{\kappa;\v}}$ and large jumps $\overline{\sL^{(\alpha)}_{\kappa;\v}}$,
$$
\sL^{(\alpha)}_{\kappa;\v} f=
\left(\int_{B_1}+\int_{B_1^c}\right)\delta^{(2)}_w f(x,\v)\kappa(t,x,\v,w)
\frac{\dif w}{|w|^{d+\alpha}}=:\widetilde{\sL^{(\alpha)}_{\kappa;\v}} f+\overline{\sL^{(\alpha)}_{\kappa;\v}} f,
$$
where 
$$
\delta^{(2)}_wf(x,\v):=f(x,\v+w)+f(x,\v-w)-2f(x,\v).
$$
By Lemma \ref{Cor72}, there is a measurable map $g(t,x,\v,w):\mR_+\times\mR^{d}\times\mR^d\times B_1\to B_1$ with 
$$
|\nabla^i_{x,\v}\nabla^j_w g(t,x,\v,w)|\leq C_{i,j}|w|^{1-j},
$$
and so that
$$
\widetilde{\sL^{(\alpha)}_{\kappa;\v}} f(x,\v)=\int_{B_1}\delta^{(2)}_{g(t,x,\v,w)}f(x,\v)\frac{\dif w}{|w|^{d+\alpha}}.
$$
Now we consider the following SDE:
$$
\dif Z_{s,t}=\int_{B_1}g(t,Z_{s,t-},w)\tilde N(\dif t,\dif w)+b(t,Z_{s,t})\dif t,\ Z_{s,s}=z.
$$
Since the coefficients are smooth and have bounded derivatives of all orders greater than $1$, it is well-known that there is a unique solution $Z_{s,t}(z)$, which
forms a $C^\infty$-stochastic flows (cf. \cite[Theorem 4.1]{Fu-Ku}). 
More precisely, it holds that for any $s<r<t$ and $z\in\mR^d$, 
$$
Z_{s,t}(z)=Z_{r,t}\circ Z_{s,r}(z),\ \  a.s.,
$$
and for any $j\in\mN$ and $p\geq 1$,
\begin{align}\label{GF9}
\sup_{z\in\mR^d}\sup_{0\leq s<t\leq T}\mE|\nabla^j_z Z_{s,t}(z)|^p<\infty.
\end{align}
Moreover, let $f:\mR_+\times\mR^{2d}\to\mR$ be a measurable function with 
$\|\nabla^j_z f\|_\infty<\infty$ for each $j\in\mN$. It is also well-known that
$$
(x,\v)\mapsto u(s,x,\v):=\int^T_s\e^{\lambda(s-t)}\mE f(t,Z_{s,t}(x,\v))\dif t\in \cap_m \sC^m(\mR^{2d})
$$
solves the following PDE:
$$
\p_s u+\widetilde\sL^{(\alpha)}_{\kappa;\v} u+b\cdot\nabla u-\lambda u=f,\ \ u(T,\cdot)=0.
$$
From the representation, by the chain rule and \eqref{GF9}, it is easy to see that for any $m\in\mN$,
there is a constant $C>0$ such that for all $s\in[0,T]$ and $\lambda\geq 0$,
\begin{align}\label{FD1}
\|u(s)\|_{\sC^m}\leq C\int^T_s\e^{\lambda(s-t)}\|f(t,\cdot)\|_{\sC^m}\dif t
\leq C\int^T_s\|f(t,\cdot)\|_{\sC^m}\dif t,
\end{align}
and by the definition of $\overline{\sL^{(\alpha)}_{\kappa;\v}}u$,
$$
\|\overline{\sL^{(\alpha)}_{\kappa;\v}} u\|_{\sC^m}\leq C\|u\|_{\sC^m}.
$$
Next we consider the following Picard's iteration: Fix $m\geq 2$. Let $u_0\equiv 0$. For $n\in\mN$,
let $u_n$ be the unique $\sC^m$-valued solution of the following PDE:
\begin{align}\label{GK3}
\p_s u_n+\widetilde{\sL^{(\alpha)}_{\kappa;\v}} u_n+b\cdot\nabla u_n-\lambda u_n=f-\overline{\sL^{(\alpha)}_{\kappa;\v}} u_{n-1},\ \ u_n(T,\cdot)=0.
\end{align}
By \eqref{FD1} we have
$$
\|u_n(s)\|_{\sC^m}\leq C\int^T_s\|f(t)-\overline{\sL^{(\alpha)}_{\kappa;\v}} u_{n-1}(t)\|_{\sC^m}\dif t
\leq C\int^T_s\Big(\|f(t)\|_{\sC^m}+\|u_{n-1}(t)\|_{\sC^m}\Big)\dif t,
$$
which yields by Gronwall's inequality that for some $C>0$,
\begin{align}\label{GK2}
\sup_n\sup_{s\in[0,T]}\|u_n(s)\|_{\sC^m}\leq C\sup_{s\in[0,T]}\|f(s)\|_{\sC^m}.
\end{align}
Similarly, we also have
$$
\|u_n(s)-u_k(s)\|_{\sC^m}\leq  C\int^T_s\Big(\|u_{n-1}(t)-u_{k-1}(t)\|_{\sC^m}\Big)\dif t,\ \ n,k\in\mN.
$$
By \eqref{GK2} and Fatou's lemma, we have
$$
\limsup_{n,k\to\infty}\sup_{t\in[s,T]}\|u_n(t)-u_k(t)\|_{\sC^m}\leq  C\int^T_s\limsup_{n,k\to\infty}\Big(\|u_{n-1}(t)-u_{k-1}(t)\|_{\sC^m}\Big)\dif t,
$$
which yields by Gronwall's inequality again
$$
\limsup_{n,k\to\infty}\sup_{t\in[0,T]}\|u_n(t)-u_k(t)\|_{\sC^m}=0.
$$
Therefore, there exists a $u\in C([0,T];\sC^m)$ so that
$$
\lim_{n\to\infty}\sup_{t\in[0,T]}\|u_n(t)-u(t)\|_{\sC^m}=0.
$$
By taking limits for equation \eqref{GK3}, we obtain the existence of $\sC^m$-valued solution to PDE \eqref{PDE0}. Since $m$ is arbitrary, we complete the proof.
\end{proof}

\end{document}